\documentclass[11pt]{article}
\usepackage{epsfig}
\usepackage{latexsym,graphicx}
\title{\bf A periodicity criterion and the section problem on the Mapping Class Group}

\author{ Patrice LE CALVEZ \thanks{ 
Institut de Math\'ematiques de Jussieu, UMR 7586 CNRS. Universit\'e Pierre et Marie Curie. Case 247, 4 Place Jussieu,
75005 Paris Cedex, France. E-mail: lecalvez@math.jussieu.fr
} }

 \date{February 14th, 2012}

\def\picture #1 by #2 (#3){
  \vbox to #2{
    \hrule width #1 height 0pt depth 0pt
    \vfill
    \special{picture #3} 
    }
  }

\def\scaledpicture #1 by #2 (#3 scaled #4){{
  \dimen0=#1 \dimen1=#2
  \divide\dimen0 by 1000 \multiply\dimen0 by #4
  \divide\dimen1 by 1000 \multiply\dimen1 by #4
  \picture \dimen0 by \dimen1 (#3 scaled #4)}
  }

\def\R{{\bf R}}
\def\Z{{\bf Z}}
\def\N{{\bf N}}
\def\Q{{\bf Q}}
\def\T{{\bf T}}
\def\A{{\bf A}}
\def\C{{\bf C}}

\def\a{\alpha}

\begin{document}
\def\eqalign#1{\null\,\vcenter{\openup\jot
\ialign{\strut\hfil$\displaystyle{##}$&
$\displaystyle{{}##}$\hfil \crcr #1\crcr }}\,}

\def\eqalignno#1{\displ@y \tabskip=\@centering
\halign to\displaywidth{\hfil$\@lign\displaystyle{##}$
\tabskip=0pt &$\@lign\displaystyle{{}##}$

\hfil\tabskip=\@centering
$\llap{$\@lign##$}\tabskip=Opt\crcr #1\crcr}}

\thispagestyle{empty}
\maketitle

\noindent{\large {\bf Abstract:}} Some years ago, V. Markovic proved that there is no section of the Mapping Class Group for a closed surface of genus $g\geq 6$ (in the case of homeomorphims) and more recently generalized this result with D. Saric to the case where $g\geq 2$. We will state a periodicity criterion and will use it to simplify some of the arguments given by
Markovic and Saric in the proof of their theorem. The periodicity criterion tells us that a homeomorphism of a connected surface must be periodic if the set of connected periodic open sets generates the topology of the surface.

\bigskip
\noindent{\bf Mathematics Subject Classification (2010)~:} \enskip
37E30, 37E45, 37B45

\bigskip
\noindent{\bf Key-words~:}\enskip Mapping class group, Dehn twist, linear hyperbolic automorphism, upper semi-continuous decomposition, shadowing lemma, rotation number, prime end.

\bigskip
\bigskip
\bigskip
\bigskip

\noindent{\large {\bf 0. Introduction}}

\bigskip
If $M$ is an orientable closed surface, denote by  ${\rm Homeo}_+(M)$ the group of orientation preserving homeomorphisms of $M$. For every $f\in {\rm Homeo}_+(M)$ write $[f]$ for the isotopy class of $f$. The set of isotopy classes of orientation preserving homeomorphisms, denoted by  ${\rm MCG}(M)$,  inherits a natural group structure such that the projection$$\eqalign{{\cal P} : {\rm Homeo}_+(M)&\to {\rm MCG}(M)\cr f&\mapsto [f]\cr}$$ is a morphism, it is called the {\it Mapping Class Group} of $M$. In the case where $M$ is the $2$-sphere, the Mapping Class Group is trivial; in the case where $M$ is the $2$-torus, it is isomorphic to ${\rm SL}(2,\Z)$. More precisely, in this last case every isotopy class contains a unique linear automorphism and such an automorphism is naturally defined by an element of ${\rm SL}(2,\Z)$. In the case where $M$ is the $2$-torus, the map that assigns to each isotopy class the automorphism that it contains, is a {\it section} of ${\cal P}$ in the following sense: it is a morphism
$${\cal E} :  {\rm MCG}(M)\to{\rm Homeo}_+(M)$$ such that
$${\cal P}\circ {\cal E}={\rm Id}_{{\rm MCG}(M)}.$$ A natural question, set by Thurston, and that can be found in Kirby's list of problems \cite{Kir} is whether such a section exists in the case where the genus of $M$ is larger than $1$. 
The same problem can be set in a smoother category. In 1987, S. Morita \cite{Mor1}, \cite{Mor2} proved that such a section does not exist for diffeomorphisms of class $C^2$, provided that $g\geq 5$. In 2007, V. Markovic \cite{Mark} proved that such a section does not exist for homeomorphisms, provided that $g\geq 6$, and extended this result to the general case $g\geq 2$  in a recent paper with D. Saric \cite{MarkS}. Summarizing, one gets:

\bigskip
\noindent{\sc Theorem A:}\enskip\enskip {\it If $M$ is an orientable closed surface of genus $g\geq 2$, there is no morphism   $${\cal E} :  {\rm MCG}(M)\to{\rm Homeo}_+(M)$$ such that $${\cal P}\circ {\cal E}={\rm Id}_{{\rm MCG}(M)}.$$ }

\bigskip
 After Markovic's paper, alternative proofs appeared. A proof by J. Franks and M. Handel \cite{FrHa1}  in the $C^1$ category working for $g\geq 3$ and a proof by S. Cantat and D. Cerveau \cite{CanCe} in the real analytic category working for $g\geq 2$. Smoothness is essential in Morita's and Cantat-Cerveau's proofs. Cantat and Cerveau deduced their theorem from the following fact: if $F$ is a real  analytic surface diffeomorphism with positive entropy, the group generated by $F$ has finite index in the real analytic centralizer of $F$. They make use of the classical result of Katok  about existence of homoclinic intersections for $C^{1+\alpha}$ diffeomorphisms. Franks-Handel's proof has a more topological flavour: it is a consequence of a topological fixed point theorem obtained by Caratheodory's prime end theory but eventually requires the differentiability because it makes use of Thurston's  Stability Theorem. By requiring only topological arguments, the proof given in  \cite {Mark} is a real  ``tour de force''. Let us conclude this brief review with two remarks.

 \smallskip
 \noindent-\enskip\enskip The proofs given in 
\cite{CanCe},  \cite{FrHa1}, \cite {Mark},   \cite{MarkS} share a common idea: looking at the centralizer of a homeomorphism (or diffeomorphism) that is isotopic to a ``local Anosov homeomorphism'' (a homeomorphism that ``coincides'' with a hyperbolic torus automorphism on a punctured torus and with the identity outside).

 \smallskip
 \noindent-\enskip\enskip In contrast with Markovic's, Markovic-Saric's and Franks-Handel's proofs, Morita's and Cantat-Cerveau's proofs extend to finite index subgroups of ${\rm MCG}(M)$ and the question is still open in the $C^0$ and $C^1$ category whether every finite index subgroup of ${\rm MCG}(M)$ admits a section.

\bigskip

Now, let us summarize the main ideas of \cite{Mark} and \cite{MarkS} (the precise definitions of the mathematical objects will be given later). The proof is by contradiction, one supposes that there exists a section ${\cal E} :  {\rm MCG}(M)\to{\rm Homeo}_+(M)$ and wants to find a contradiction.

Let us begin with some simple observations. If $F$ and $G$ are orientation preserving homeomorphisms with disjoint supports, then $F$ and $G$ commute, which implies that the classes $[F]$ and $[G]$ commute and therefore that the homeomorphisms ${\cal E}([F]) $ and ${\cal E}([G])$ also commute. Note also that if $F$ has order $q$, then ${\cal E}([F]) $ has order $q$.

Dehn twists play an essential role in the proofs.  If $\beta$ is a simple closed curve non homotopic to zero, one can define a Dehn twist $F^*_{\beta}$ supported  in an annular neighborhood $A_{\beta}$ of $\beta$. The isotopy class $[F^*_{\beta}]$ depends only on the free homotopic class of $\beta$ and one can define $F_{\beta}={\cal E}([F^*_{\beta}])$. If $\beta$ and $\beta'$ are two disjoint simple closed curves, the annular neighborhoods can be chosen to be disjoint, which implies that $F_{\beta}$ and $F_{\beta'}$ commute. Suppose now that $\beta$ and $\beta'$ have a unique point of intersection and that the intersection is transverse. In that case, the {\it Artin relations} can be written:
$$[F^*_{\beta}][F^*_{\beta'}][F^*_{\beta}]=[F^*_{\beta'}][F^*_{\beta}][F^*_{\beta'}]$$which implies that
$$F_{\beta}\circ F_{\beta'}\circ F_{\beta}=F_{\beta'}\circ F_{\beta}\circ F_{\beta'}.$$
The fundamental idea of Markovic is to show that the dynamics of each $F_{\beta}$ is so ``close'' to the dynamics of $F^*_{\beta}$ that the previous relations cannot be satisfied, at least for every couple $(F_{\beta}, F_{\beta'})$ (we will conclude this introduction by proving that there is no section that sends the isotopy classes of Dehn twists onto Dehn twists). The word ``close'' in the previous sentence means that there exists an invariant open set $U_{\beta}$ isotopic to $M\setminus A_{\beta}$ such that the restriction of $F_{\beta}$ to $U_{\beta}$ is a {\it cellular extension} of the identity map on $U_{\beta}$. This means that it is a topological extension of the identity with cellular fibers (intersection of nested sequences of topological closed disks). By  a classical result of Moore \cite{Moo}, it is equivalent to say that there exists a semi-continuous decomposition of $U_{\beta}$ by invariant cellular sets. Moreover, the set $U_{\beta}$ and the decomposition are canonical in the following sense: for every homeomorphism $F\in{\rm Im} ({\cal E})$, one has $F(U_{\beta})=U_{F(\beta)}$ and $F$ sends fibers of the first extension onto fibers of the second extension. The contradiction obtained by Markovic in the case where $g\geq 6$ makes use of the many symmetries that occur, the contradiction obtained by Markovic and Saric in the  general case arises from a different relation that we will explain now in the case where $g$ is even (the odd case being slightly different). 

\medskip

One can find a family $(\beta_j)_{-g\leq j\leq g}$ of simple closed curves such that:

\smallskip
\noindent-\enskip\enskip $\beta_0$ separates $M$ in two homeomorphic one punctured surfaces of genus $g/2$;

\smallskip
\noindent-\enskip\enskip the $\beta_j$, $j<0$, are included in one of the connected component of $M\setminus\beta_0$, and the $\beta_j$, $j>0$, are included in the other one;

\smallskip
\noindent-\enskip\enskip there exists an involution  $I^*\in {\rm Homeo}_+(M)$ that fixes $\beta_0$, that permutes the two connected components of $M$, and that permutes $\beta_j$ and $\beta_{-j}$;

\smallskip
\noindent-\enskip\enskip if $\vert j-j'\vert >1$, then $\beta_j$ are $\beta_{j'}$ disjoint;

\smallskip
\noindent-\enskip\enskip  if $j$, $j'$ are different from $0$ and if $\vert j-j'\vert =1$, then $\beta_j$ and $\beta_{j'}$ have a unique point of intersection and the intersection is transverse.

\medskip One deduces that

\smallskip
\noindent-\enskip\enskip $I={\cal E}(I^*)$ is an involution;

\smallskip
\noindent-\enskip\enskip $I$ conjugates $F_{\beta_j}$ to $F_{\beta_{-j}}$;

\smallskip
\noindent-\enskip\enskip $F_{\beta_0}$ commutes with every  $F_{\beta_j}$;

\smallskip
\noindent-\enskip\enskip $F_{\beta_j}$ and  $F_{\beta_j}$ commute if $\vert j-j'\vert>1$ . 

\smallskip
\noindent The contradiction will follow from the equality
$$( F_{\beta_1}\circ \dots \circ  F_{\beta_g})^{2g+2}=F_{\beta_0}=(F_{\beta_{-1}}\circ \dots \circ F_{\beta_{-g}})^{2g+2},$$
which is a consequence of the following Artin type relation  (see Farb-MargalitÊ \cite{FaMarg}):
$$([ F^*_{\beta_1}]\dots  [ F^*_{\beta_g}])^{2g+2}=[F^*_{\beta_0}]=([F^*_{\beta_{-1}}] \dots [F^*_{\beta_{-g}}])^{2g+2}.$$

Let us explain now how Markovic constructs this set $U_{\beta}$ on which $F_{\beta}$ is a cellular extension of the identity. Let us begin with a simpler situation. Consider an orientation preserving hyperbolic automorphism $L$ of the torus $\T^2=\R^2/\Z^2$ with exactly one fixed point. It is a classical fact, consequence of the well known Shadowing Lemma (see Franks \cite{Fr}),  that the centralizer of $L$ in ${\rm Homeo}_+(\T^2)$ is the cyclic group generated by $L$. In particular, every homeomorphism $F$ isotopic to the identity that commutes with $L$ is trivial. Another classical consequence of the Shadowing Lemma,  is the fact that every homeomorphism $G$ of $\T^2$ isotopic to $L$ is an extension of $L$ and that the semi-conjucacy $H$ is uniquely defined if homotopic to the identity (see Franks \cite{Fr}). It is not difficult to show that the fibers of $H$ are {\it acyclic} (which means that their connected components are cellular). Suppose now that $F$ is a homeomorphism isotopic to the identity that commutes with $G$. By uniqueness of $H$, it is easy to see that each fiber is invariant by $F$. What is much more surprising and is proven in \cite{Mark} is that every connected component of a fiber is also invariant. The decomposition in connected components of fibers of $H$ is an upper semi-continuous decomposition in cellular sets, which implies that $F$ is a cellular extension of the identity. Such results can be generalized to the case where $L$ is a local Anosov homeomorphism of $M$. More precisely, suppose that $M=N\cup A\cup T$ can be decomposed in three invariant surfaces with boundary (with disjoint interiors) such that:

\smallskip
\noindent-\enskip\enskip $T$ is a one punctured compact torus and $L\vert_T$ is conjugate to an orientation preserving hyperbolic automorphism of $\T^2$ with one fixed point, when one blows up this fixed point;

\smallskip
\noindent-\enskip\enskip $N$ is a compact one punctured surface of genus $g-1$ and $L$ fixes every point of $N$;

\smallskip
\noindent-\enskip\enskip $A$ is an annulus and $L\vert_A$ is isotopic to the identity relative to the set of fixed points that lie on the boundary.

\smallskip
\noindent Suppose now that $G$ is isotopic to $L$ and that $F$ commutes with $G$ and is isotopic to a homeomorphism that fixes every point of $T$. Then there exists an invariant open set $U$ which ``contains'' the homology of $T$ and an upper semi-continuous decomposition of $U$ by $F$-invariant acyclic sets.

It remains to observe that for every simple closed curve $\beta$ non homotopic to zero and for every compact one punctured torus $T\subset M\setminus A_{\beta}$, there exists a local Anosov homeomorphism $L^*$ supported on $T$. As it commutes with  $F^*_{\beta}$, one deduces that ${\cal E}(L^*)$ commutes with $F_{\beta}$. Consequently, there exists an  open set $U$ which contains the homology of $T$ and an upper semi-continuous decomposition of $U$ by $F_{\beta}$-invariant acyclic sets. The choice of another torus $T'$ gives us another open set $U'$. It remains to prove that the ``intersection'' of two such decompositions defined respectively on $U$ and $U'$ defines an upper semi-continuous decomposition by invariant acyclic sets on $U\cup U'$, that one can construct in that way an upper semi-continuous decomposition by invariant acyclic sets on an open set isotopic to $M\setminus A_{\beta}$, to prove that such a decomposition can be constructed to be canonical and finally to prove the surprising fact that the decomposition induced by taking connected components is also invariant.

\bigskip
The two articles \cite{Mark} and \cite{MarkS} are not easy to read. The ideas are wonderful but the proofs are tough and there is sometimes a lack of precision both in the statements and in the proofs (for example Lemma 2.1 of \cite{Mark} is untrue stated at it is, which does not help the reader to understand the arguments of Markovic in his study of upper semi-continuous decompositions of surfaces). The goal of this work was to re-write the proof of both papers, to get a more precise version. Indeed, as it is said above, one does not know if the lifting problem can be generalized to finite index subgroups and we hope that  our effort could help. Of course, there would have been no interest  in following too closely the original articles and we have tried to give more conceptual arguments (at least we hope!). It appears that this can be done by using a  periodicity theorem (Theorem B below) that does not appear in \cite{Mark} and \cite{MarkS}. To write an article as self-contained as possible we have added some proofs of more classical results.

\medskip

Let us explain now the plan of the article.

\medskip
In Section $1$, we will state the following periodicity criterion, whose proof needs Brouwer's Lemma on Translation Arcs (a result about homeomorphisms of the plane). A natural question is whether this result can be extended to higher dimension.
 
\bigskip
\noindent{\sc Theorem B:}\enskip {\it  A homeomorphism of a connected surface $M$  is periodic if and only if the set of periodic connected open sets generates the topology of $M$.}

\bigskip
In section 2 we recall the definition of an upper semi-continuous decomposition. After stating classical results about such decompositions in the most general framework, we will introduce  {\it local upper semi-continuous decompositions} in the case of a locally compact Hausdorff space, with their  {\it domain}. We will also define the {\it intersection} of local upper semi-continuous decompositions. After that, we will focus on decompositions of surfaces and will introduce {\it essential, acyclic,} and {\it cellular decompositions}. Cellular decompositions are local decompositions whose elements are cellular sets; acyclic decompositions are local decompositions whose elements have cellular connected components (we do not follow the terminology given in \cite{Mark}); essential decompositions are local decompositions such that no connected component of the domain is included in a disk of $M$. Moreover if $F$ is a homeomorphism of a surface, we will introduce {\it invariant} decompositions which are local upper semi-continuous decompositions whose elements are invariant by $F$. A fact, easy to prove but fundamental, is that the intersection of essential acyclic decompositions is still acyclic (in opposition to non essential acyclic decompositions !). This permits us to define the {\it canonical invariant essential acyclic decomposition} of a homeomorphism which is the intersection of all the invariant essential cellular decompositions and which is the finest among all such decompositions. A priori, there is no reason why something similar could be done for cellular decompositions because the intersection of a family of cellular sets is acyclic but not necessarily cellular. The surprising fact that such a canonical cellular decomposition exists is one of the key result of \cite{Mark} (but it is not stated in the same words).  Using Theorem B, we will give a prove of its existence:

\bigskip
\noindent{\sc Theorem C:}\enskip\enskip {\it Let $F$ be a homeomorphism of an oriented (not necessarily closed) surface $M$.  The intersection of all the invariant essential cellular decompositions is itself an invariant essential cellular decomposition, it is the finest among all such decompositions, we call it the {\rm canonical invariant essential cellular decomposition} of $F$.}

\bigskip

For a given homeomorphism, the domain of the invariant essential acyclic (or cellular) decomposition is usually empty. However, as it is explained in \cite{Mark}, its domain can be large if the centralizer is wide. The main result of Section 3 gives a precise meaning to this. The theorem below  is proved in \cite{Mark}. Instead of the topological arguments that are used in \cite{Mark} we will apply in a more classical way the  {\it Shadowing Lemma} for hyperbolic automorphisms of $\T^2$.

\bigskip
\noindent{\sc Theorem D:}\enskip\enskip {\it Let $M$ be an orientable  closed surface with no boundary, $T\subset M$ a compact one punctured torus, $\check T$ the $2$-torus obtained from $M$ by identifying $\overline {M\setminus T}$ to a point and $P~: M\to \check T$ the natural projection. Let $F$ be a homeomorphism of $M$ such that:

\smallskip 
\noindent{\bf i)}\enskip\enskip $F$ is isotopic to a homeomorphism that fixes every point of $T$;

\smallskip 
\noindent{\bf ii)}\enskip\enskip  there exists a homeomorphism that commutes with $F$ and is isotopic in $M$ to a local Anosov homeomorphism supported on $T$.

\smallskip
Then, there exists a connected component $U$ of the domain of the canonical invariant essential acyclic decomposition of $F$ such that the map $P^*: H^1(\check T,\Z)\to H^1(U,\Z)$ is injective.}

\bigskip

Finally, in section 4 we will give the proof of Theorem A. The proof is close in spirit to what is done in \cite{MarkS} but is quite different and (we hope) simpler (for example there is no need to introduce pseudo-Anosov maps). The final argument is the following result (easy to prove) on homeomorphims of the circle:

\bigskip
\noindent{\sc Proposition E:}\enskip\enskip {\it There is no family $(\Phi_j)_{n\leq j\leq m}$, $n<0<m$, of orientation preserving homeomorphisms of the circle $S^1$ such that:

 \smallskip
 \noindent-\enskip each $\Phi_j$, $j<0$, commutes with each $\Phi_{j'}$, $j'>0$;
  
  \smallskip
 \noindent-\enskip $\Phi_j$ has a fixed point if and only if $j\not=0$;

  \smallskip
 \noindent-\enskip one has $\Phi_1\circ \dots \circ  \Phi_m=\Phi_0=\Phi_{-1}\circ \dots \circ \Phi_{-n}$.}
 
\bigskip This proposition implies immediatly that there is no section $\cal E$ of $\cal P$ that sends the isotopy classes of Dehn twists onto Dehn twists. Indeed,  suppose that such a section exists and look at the case where  the genus is even. We keep the notations $I$ and $F_{\beta_j}$, $-g\leq j\leq g$, given at the beginning of the introduction. The set $\Gamma$ of periodic points of period $2$ of $F_{\beta_0}$ is a circle invariant by each $F_{\beta_j}$, $-g\leq j\leq g$, because $F_{\beta_j}$ commutes with $F_{\beta_0}$. This circle $\Gamma$, being homotopic to $\beta_0$, is not homotopic to $\beta_j$, if $j\not=0$, and therefore is not contained in the support of $F_{\beta_j}$. Consequently, it contains a fixed point of $F_{\beta_j}$. Proposition E, applied to the restrictions of the $F_{\beta_j}$to $\Gamma$  tells us that the condition 
 $$( F_{\beta_1}\circ \dots \circ  F_{\beta_g})^{2g+2}=F_{\beta_0}=(F_{\beta_{-1}}\circ \dots \circ F_{\beta_{-g}})^{2g+2}$$
cannot be satisfied.  To use this argument in the general case, one must prove that the group generated by the $F_{\beta_j}$, $-g\leq j\leq g$, acts on a circle $\Gamma$ in such a way that $F_{\beta_j}$ has a fixed point if and only if $j\not=0$. This circle will be a {\it prime end circle}. More precisely we will construct a compact connected set $X$ that satisfies the following properties:

 \smallskip
 \noindent-\enskip $X$ is invariant by $I$ and by each $F_{\beta_j}$, $-g\leq j\leq g$;

 \smallskip
 \noindent-\enskip the complement of $X$ is the union of two one punctured surfaces $V_-$, $V_+$ of genus $g/2$ that are permuted by $I$ and invariant by each $F_{\beta_j}$, $-g\leq j\leq g$;

 \smallskip
 \noindent-\enskip If $V_-\sqcup S^1$ is the {\em prime end compactification} of $V_-$, then the rotation number of each map $F_{\beta_i}\vert_{V_-}$, $j\not=0$ induced on the circle $S^1$ is $0+\Z$ while the rotation number of $F_{\beta_0}\vert_{V_-}$ is $1/2+\Z$.
 
 \bigskip
  \bigskip
In all the article, a {\it surface} will be a Hausdorff topological surface without boundary; an {\it open disk} of a surface a set homeomorphic to  ${\bf D}=\{ z\in\C\, \,\vert \,\vert z\vert <1\}$; a {\it closed disk} a set homeomorphic to $\overline{\bf D}=\{ z\in\C\, \vert\, \vert z\vert \leq1\}$, an {\it open annulus} a set homeomorphic to $\R/\Z\times\R$, a {\it closed annulus} a set homeomorphic to $\R/\Z\times[0,1]$. If $Y$ is a subset of a topological space $X$, we will denote by $\overline X$, ${\rm Int}(X)$ and ${\rm Fr}(X)$ the closure, the interior and the frontier respectively.

\bigskip
Trying to understand the article of Markovic was the goal of a weekly workshop that took place in Paris in 2007-2008 as part of the project  Symplexe (ANR-06-BLAN-0030-01). I would like to thank the participants of this workshop and particulary Fran\c cois B\'eguin, Sorin Dumitrescu and Fr\'ed\'eric Le Roux for the fruitful discussions we had at that time.

\bigskip
\bigskip

\noindent{\large {\bf 1. A periodicity criterion}}

\bigskip
The main goal of this first section is to prove the following result~:

\bigskip
\noindent{\sc Theorem 1.1:}\enskip {\it  Let $F$ be a homeomorphism of a connected surface $M$. Then the following assertions are equivalent:

\medskip
\noindent-\enskip\enskip the topology of $M$ is generated by periodic  connected open sets {\bf (P)};

\medskip
\noindent-\enskip\enskip there exists $q\geq 1$ such that $F^{q}={\rm Id}_M$ {\bf (Q)}.}

\bigskip

The idea of the proof is to start with a homeomorphism $F$ of a connected surface $M$, that satisfies {\bf (P)}. Using finitely many elementary modifications (we will define the meaning later) we will construct an orientation preserving  homeomorphism $F'$ of the $2$-sphere $S^2$ that has at least three fixed points and that satisfies {\bf (P)}. The main result of this section (Proposition 1.5) is that such a homeomorphism $F'$ must coincide with the identity on $S^2$. By definition of elementary modifications, this will immediatly imply that $F$ is periodic.

\medskip
Let us begin with simple observations. Let $F$ be a homeomorphism of a surface $M$ that satisfies {\bf (P)}, then:

\smallskip
\noindent-\enskip\enskip every iterate $F^k$ of $F$, $k\in\Z$, satisfies {\bf (P)};

\smallskip
\noindent-\enskip\enskip the restriction of  $F$ to any invariant open set satisfies {\bf (P)};

\smallskip
\noindent-\enskip\enskip every lift of $F$ to a finite covering space of $M$ satisfies {\bf (P)}.

\medskip 
Denote by ${\cal V}_F$ the set of periodic  connected open sets  that are included in a closed disk of $M$. 
The following remarks are obvious:

\smallskip
\noindent-\enskip\enskip each set $V\in{\cal V}_F$ is relatively compact;

\smallskip
\noindent-\enskip\enskip the set ${\cal V}_F$ generates the topology of $M$;
 
\smallskip
\noindent-\enskip\enskip the equality ${\cal V}_F={\cal V}_{F^k}$ is  true for every $k\geq 1$;

\smallskip
\noindent-\enskip\enskip if $\widetilde M$ is a covering space of $M$, then, for every $V\in{\cal V}_F$, the covering projection  $\Pi~:\widetilde M\to M$ induces a homeomorphism between every connected component $\widetilde V$ of $\Pi^{-1}(V)$ and $V$;

\smallskip
\noindent-\enskip   the set $\widetilde{{\cal V}_F}$ of all open sets $\widetilde V$ constructed as above generates the topology of $\widetilde M$.

\bigskip
\bigskip

\noindent{\bf Preliminary results}

\bigskip

\noindent{\sc Proposition 1.2~:}\enskip{\it Let $M$ be a connected surface, $F$ a homeomorphism of $M$ that satisfies {\bf (P)} and $\widetilde F$ a lift of $F$ to a connected covering space $\widetilde M$ of $M$. Then, the following alternative holds:

\smallskip
\noindent-\enskip\enskip every set
$\widetilde V\in\widetilde{{\cal V}_F}$ is $\widetilde F$-periodic, in this case  $\widetilde F$ satisfies  {\bf (P)} and one has $\widetilde{{\cal V}_F}\subset{\cal V}_{\widetilde F}$;

\smallskip
\noindent-\enskip\enskip none of the sets $\widetilde V\in\widetilde{{\cal V}_F}$ is $\widetilde F$-periodic, in that case  $\widetilde F$ has no periodic point.
}

\bigskip
\noindent{\it Proof.}\enskip Denote by $\widetilde W_p$ the union of the sets $\widetilde V\in\widetilde{{\cal V}_F}$ that are periodic and by $\widetilde W_n$ the union of the sets $\widetilde V\in\widetilde{{\cal V}_F}$ that are not periodic. We get two open sets and we want to prove that one of them is empty. The space $\widetilde M$ being connected, it is sufficient to prove that $\widetilde W_p$ does not intersect $\widetilde W_n$. Equivalently, one must prove that every set $\widetilde V'\in\widetilde{{\cal V}_F}$ that meets a periodic set $\widetilde V\in\widetilde{{\cal V}_F}$ is itself periodic. Denote by $V'$ the image of $\widetilde V'$ in $M$ by the covering projection. One can find an integer $k\geq 1$ such that $\widetilde V$ is fixed by $\widetilde F^k$ and $V'$ fixed by $F^k$. If ${\cal T}$ is the deck transformation group, there are finitely many elements in the ${\cal T}$-orbit of $\widetilde V'$ that meet $\widetilde V$, because $\widetilde V$ and $\widetilde V'$ are relatively compact. The homeomorphism  $\widetilde F^k$ inducing a permutation on this finite set, one deduces that $\widetilde V'$ is periodic. \hfill$\Box$

\bigskip
\noindent{\sc Proposition 1.3~:}\enskip {\it  Let $F$ be a homeomorphism of a surface $M$. If $F$ satisfies  {\bf (P)}, then the set of periodic points is dense.}

\bigskip
\noindent{\it Proof.}\enskip An immediate consequence of the hypothesis is the fact that the natural action of $F$ on the set of connected components of $M$ has finite orbits. Thus, one can suppose that $M$ is connected. Replacing $M$ with the orientable two-folds covering space in the case where $M$ is non orientable, one can suppose that $M$ is orientable. Replacing $F$ with $F^2$ in the case where $F$ reverses the orientation, one can suppose that $F$ preserves the orientation. We want to prove that every open disk $W\subset M$ contains a periodic point. Let us choose $V\in{\cal V}_F$ such that $\overline V\subset W$, and $k\geq 1$ such that $F^k(V)=V$. The connected component $X$ of $M\setminus V$ that contains $M\setminus W$ is $F^k$-periodic. Indeed $F^k$ acts naturally on the set of connected components of $M\setminus V$ and the $F^k$-orbit of $X$ is finite because $X$ has no empty interior and $F^k$ satisfies {\bf (P)}. Replacing $F$ with an iterate $F^{kk'}$, $k'\geq 1$, one can suppose that both $V$ and $X$ are invariant by $F$. The set $V'=M\setminus X$, which is invariant by $F$, is the union of $V$ and of the compact connected components of $W\setminus V$: it is an open disk. To end the proof,  we will prove that $F\vert_{V'}$ has  a fixed point. The argument, originally due to Brown \cite{Brow}, is frequently used in the dynamical study of homeomorphisms of surfaces.  If $F$ fixes every point of $V'$, we are done, otherwise one can find a periodic connected open set $V''\subset V'$ such that $F(V'')\cap V''=\emptyset$. Denote by $q$ the smallest positive integer such that $f^q(V'')\cap V''\not=\emptyset$ and fix $x\in V''\cap F^{-q}(V'')$.  One can find a homeomorphism $h$ of $V'$ supported on $V''$ that sends $F^q(x)$ on $x$. The point $x$ is a periodic point of $h\circ F\vert_{V'}$, of period $q$. The well-known Brouwer's Lemma on Translation Arcs \cite{Brou} tells us that $h\circ F\vert_{V'}$ has a fixed point in $V'$. Observe now that this point is fixed by $F$ because $V''$ does not meet its image.  \hfill$\Box$

\bigskip
\noindent{\bf Remark~:}\enskip\enskip Instead of Brouwer's Lemma, one can use Cart\-wright-Little\-wood's Fixed Point Theorem \cite{CarL}. Let us explain why, by keeping the notations of the previous proof.
Replacing $F$ by an iterate $F^m$ if necessary,  one can suppose that the connected component $U$ of $M\setminus\overline V$ that contains $M\setminus W$ is fixed by $F$. The set $Y=M\setminus U$ is compact, connected, included in $W$ and invariant by $F$. Moreover $W\setminus Y$ is connected. This implies that $Y$ is cellular. Cartwright-Littlewood's Fixed Point Theorem asserts that $F$ has a fixed point inside $Y$.

\bigskip
\bigskip

\noindent{\bf Elementary modifications}

\bigskip 

Let us introduce  a last definition. Let $F$ and $F'$ be homeomorphisms of  connected surfaces $M$ and $M'$ respectively. Say that $F'$ is an {\it  elementary modification} of $F$ in any of the following cases:

\smallskip
\noindent-\enskip\enskip $F'=F^k$ is a positive iterate of $F$;

\smallskip
\noindent-\enskip\enskip $F'$ is a lift of $F$ to a covering space of $M$;

\smallskip
\noindent-\enskip\enskip $F'$ is the restriction of $F$ to the complement of a periodic orbit of $F$;

\smallskip
\noindent-\enskip\enskip $F$ is the restriction of $F'$ to the complement of a periodic orbit of $F'$.

\bigskip
\noindent{\sc Proposition 1.4~:}\enskip {\it  Let $F$ and $F'$ be homeomorphisms of  connected surfaces $M$ and $M'$ respectively.

\smallskip
\noindent{\bf i)}\enskip\enskip If $F'$ is an elementary modification of $F$ and if $F$ satisfies {\bf (P)}, then $F'$ satisfies {\bf (P)} except if $F'$ is a periodic point free lift of $F$ to a covering space of $M$.

\smallskip
\noindent{\bf ii)}\enskip\enskip If $F'$ is obtained from $F$ after finitely many elementary modifications and if $F'$ is periodic, then $F$ itself is periodic.

}

\bigskip
\noindent{\it Proof.}\enskip The assertion {\bf ii)} is clear. If we want to prove {\bf i)}, the only case to consider which is not obvious, is the case where $F$ is the restriction of $F'$ to the complement of a periodic orbit $O$ of $F'$. One must prove that every point $z\in O$ has a fundamental system of $F'$-periodic connected open neighborhoods. Denote by $q$ the period of $O$ and fix a neighborhood $W$ of $z$ that does not contain any other point of $O$. One must find a $F'$-periodic connected open neighborhood $V$ of $z$ that is included in $W$. Consider a closed disk $D\subset W$ whose interior contains $z$. By hypothesis, the boundary $\partial D$ can be covered by a family $(V_i)_{i\in I}$ of  periodic connected open sets included in $M\setminus  O$ and relatively compact in $M\setminus  O$. By compactness of $\partial D$, one can suppose that $I$ is finite. The set $K=\bigcup_{i\in I, k\in\Z} F{}^{k}(\overline V_i)$ is compact because there are finitely many $F{}^{k}(\overline V_i)$ in the union and each of these sets is compact. More precisely, $K$ is an invariant compact set included in $M\setminus O$ and containing $\partial D$. The connected component $V$ of $M\setminus K$ that contains $z$ is included in $W$ and does not contain any other point of $O$ because 
$\partial D$ is included in $K$. One deduces that $V$ is $F'$-periodic, of period $q$. \hfill$\Box$

\bigskip
\bigskip

\noindent{\bf Proof of the theorem}

\bigskip

\noindent{\sc Proposition 1.5~:}\enskip {\it  Let $F$ be an orientation preserving homeomorphism of the sphere $S^2$ that satisfies {\bf (P)} and that has at least three fixed points. Then $F$ is the identity.
}

\bigskip
\noindent{\it Proof.}\enskip We will give a proof by contradiction and suppose that $S^2\setminus{\rm Fix}(F)$ is not empty. We fix a connected component $U$ of $S^2\setminus{\rm Fix}(F)$. A result of Brown and Kister \cite{BrowKis} asserts that $U$ is invariant by $F$, which implies that we can find in $U$ a path $\gamma $ that joins a point $z$ to its image $F(z)$.  As the restriction $F\vert_U$ satisfies {\bf (P)}, one can cover  $\gamma$ by a family $(V_i)_{i\in I}$ of  periodic connected open sets included in $U$ and relatively compact in $U$. Here again, one can suppose that $I$ is finite and minimal: there is no subfamily that covers $\gamma$. The set $W=\bigcup_{i\in I, k\in\Z} F{}^{k}(V_i)$ is an invariant open set. It is relatively compact because each $F^k(V_i)$ is relatively compact and there are finitely many such sets in the union. It is connected because 
 each $F^k(V_i)$ is connected and meets $\bigcup_{k\in\Z} F^k(\gamma)$ which is connected and included in $W$. The fact that ${\rm Fix}(F)$ is disjoint from $\overline W=\bigcup_{i\in I, k\in\Z} F^k(\overline{V_i})$ implies that there are finitely many connected components of $S^2\setminus \overline W$ that contain a fixed point. Indeed ${\rm Fix}(F)$ is compact and covered by the connected components of $S^2\setminus \overline W$. Obviously, this implies that there are finitely many connected component of $S^2\setminus W$ that contain a fixed point. We will denote by 
$\cal X$ the set of such components. One may observe that every connected component of $S^2\setminus W$ is cellular, because $W$ is connected, and that every such component which is fixed by $F$ contains a fixed point according to Cartwright-Littlewood's Fixed Point Theorem \cite{CarL}. Therefore, $\cal X$ is the union of the connected components of $S^2\setminus W$ that are fixed by $F$ (we will not use this fact in the proof of the proposition but it helps to understand the situation).

\medskip 
The complement of $\bigcup_{X\in {\cal X}} X$ is a finite punctured sphere, invariant by $F$ and fixed point free. Indeed, every fixed point of $F$ belongs to a set $X\in{\cal X}$ by definition of ${\cal X}$. Proposition 1.3 tells us that the complement of $\bigcup_{X\in {\cal X}} X$ contains periodic  points. One deduces, using Brouwer's Lemma on Translation Arcs, that the complement of $\bigcup_{X\in {\cal X}} X$ is not a disk: there are at least two punctures. 

\medskip

Let us begin by the case where there are more than two punctures and explain why there is a contradiction. One can find a set ${\cal X}'\subset{\cal X}$ such that  $F\vert_{S^2\setminus \bigcup_{X\in{\cal X}'}X}$ is  isotopic to the identity and which is maximal relative to this property. As one knows that every orientation preserving homeomorphism of a three punctured sphere is isotopic of the identity, one deduces that  ${\cal X}'$ has at least three elements. Thus, the fundamental group of $N= S^2\setminus \bigcup_{X\in{\cal X}'} X$ is the free group with $\sharp {\cal X}'-1$ generators and its center is trivial. This implies that there exists a unique lift $\widetilde G$ of $G=F\vert_N$ to the universal covering space $\widetilde N$ that commutes with the deck transformations. The maximal property supposed on ${\cal X}'$ results in the fact that every set $X\in{\cal X}\setminus {\cal X}'$ is lifted to subsets of $\widetilde N$ that are not fixed by $F$. In other terms, $\widetilde G$ is fixed point free. According to Brouwer's Lemma on Translation Arcs, one deduces that $\widetilde G$ is periodic point free. As a consequence of Proposition 1.3 and Proposition 1.2, one deduces first that  $\widetilde G$ does not satisfy {\bf (P)} and then that there is no  $\widetilde V\in {{\cal V}_{G}}$ which is $\widetilde G$-periodic. Denote by $\cal T$ the deck transformation group and fix $\widetilde V\in\widetilde {{\cal V}_{G}}$. By definition, $\widetilde V$ lifts a $F$-periodic connected open set. Thus, there exists $\tau\in{\cal T}\setminus\{{\rm Id}_{\widetilde N}\}$ and $k\geq 1$ such that $\widetilde G^k(\widetilde  V)=\tau(\widetilde  V)$. According to Proposition 1.2, one deduces that the lift $\widetilde G^k\circ {T^{-1}}$ of $G^k$ satisfies {\bf (P)}. Therefore, for every $\tau'\in{\cal T}$, there exists $k'\geq 1$ such that $ \tau'(\widetilde V)=(\widetilde G^k\circ \tau^{-1})^{k'}(\tau'(\widetilde V))$. The homeomorphism $\widetilde G$ commutes with the deck transformations, so we have :
$$\tau'(V)=\tau^{-k'}\circ \tau'\circ \widetilde G^{kk'}(V)=\tau^{-k'}\circ \tau' \circ \tau^{k'}(V),$$ which means that $\tau^{k'}\circ \tau'=\tau'\circ \tau^{k'}$. Recall that $\cal T$ is the free group with $\sharp {\cal X}'-1$ generators. One deduces that there exists a power of $\tau'$ that commutes with every $\tau\in{\cal T}$, which is impossible.

 \medskip
It remains to study the case where ${\cal X}'$ has two elements. We write ${\cal X}'=\{X_1,X_2\}$. In this case, $N=S^2\setminus(X_1\cup X_2)$ is an annulus invariant by $F$ that does not contain fixed points of $F$. Let us choose a fixed point $z_1\in X_1$ and a fixed point $z_2\in X_2$. Set $M'=S^2\setminus\{z_1,z_2\}$. By hypothesis, there exists a third fixed point $z_3$. The lift of  $F\vert_{M'}$ to the universal covering space $\widetilde {M'}$ that fixes the lifts of $z_3$ satisfies {\bf (P)}, according to Proposition 1.2. Denote by $\Pi~:\widetilde {M'}\to M'$ the covering projection. The annulus $N$ being essential in $M'$, the set $\widetilde N=\Pi^{-1}(N)$ is the universal covering space of $N$, it is a topological plane. The restriction $\widetilde {F'}\vert_{\widetilde N}$ of $\widetilde {F'}$ to $\widetilde N$ satisfies {\bf (P)} because $\widetilde {F'}$ does. One deduces that $\widetilde {F'}\vert_{\widetilde N}$ is not periodic point free by Proposition 1.3. But $\widetilde {F'}\vert_{\widetilde N}$ is a lift of  $F\vert_N$ and $F\vert_N$ is fixed point free. This implies that $\widetilde {F'}\vert_{\widetilde N}$ is fixed point free. The contradiction comes from Brouwer's Lemma on Translation Arcs. \hfill$\Box$
\bigskip 

\noindent{\it Proof of Theorem 1.1.}\enskip Let $M$ be a connected surface and $F$ a homeomorphism of $M$ that satisfies {\bf (P)}. Using an elementary modification if necessary (replacing $M$ with an orientable two folds covering space), one can suppose that $M$ is orientable. Using another elementary modification if necessary (replacing $F$ with $F^2$), one can suppose that $F$ preserves the orientation.

In the case where $M$ is the $2$-sphere, one can use Proposition 1.3 and find an iterate $F^k$ of $F$ that has at least three fixed points. This map has been obtained after at most three elementary modifications and Proposition 1.5 tells us that it coincides with the identity.  As a consequence, one deduces that our original map $F$ is periodic.

In the case where $M$ is not the $2$-sphere, one needs three more elementary modifications. A new elementary modification (replacing $F$ with an iterate $F^k$) permits us to suppose that $F$ has a fixed point denoted by $z$. Let us write $\widetilde M$ for the universal covering space of $M$ and fix a lift $\widetilde z\in\widetilde M$ of $z$. There exists a lift $\widetilde F$ of $F$ that fixes $\widetilde z$. This lift satisfies {\bf (P)} by Proposition 1.2. As we suppose that $M$ is not the $2$-sphere, we know that $\widetilde M$ is topologically a plane. The extension of $\widetilde F$ to the one point compactification of $\widetilde M$ gives us a homeomorphism of the $2$-sphere and we know that it satisfies {\bf (P)} because of Proposition 1.4. Taking a certain power of this extension, one gets a homeomorphism of the $2$-sphere that satisfies {\bf (P)}, that has at least three fixed points, and that is obtained from our original transformation after at most six elementary modifications.  \hfill$\Box$

\bigskip
\bigskip

\noindent{\bf Remarks about periodic homeomorphisms of surfaces}

\bigskip
Proposition 1.5 applied to periodic homeomorphims permits to get a proof of very classical results about such homeomorphisms. As we will use these results later, we recall them.

\bigskip

\noindent{\sc Proposition 1.6~:}\enskip {\it  {\bf i)} The fixed point set of a non trivial periodic orientation preserving homeomorphism of a connected orientable surface is discrete.

 \smallskip
 {\bf ii)} Two non trivial periodic orientation preserving homeomorphims  of a closed connected orientable surface $M$ of genus $g\geq 2$ that are isotopic and have the same period, have the same (finite) number of fixed points.} 

\bigskip

\noindent{\it Proof.} \enskip Let us begin with the proof of {\bf i)}. Let $F$ be a periodic orientation preserving homeomorphism of a connected orientable surface $M$ such that its fixed point set has an accumulation point $z$.  If $M$ is the $2$-sphere, Proposition 1.5 tells us that $F$ is trivial, because $F$ satisfies {\bf (P)}. If $M$ is not the $2$-sphere, consider the universal covering space $\widetilde M$ of $M$ and a lift $\widetilde z\in\widetilde M$ of $z$. If $\widetilde F$ is the lift of $F$ that fixes $\widetilde z$, then $\widetilde z$ is an accumulation point of the set of fixed points of $\widetilde F$. The extension of $\widetilde F$ to the one point compactification of $\widetilde M$ is a homeomorphism of the $2$-sphere that satisfies {\bf (P)} and has infinitely many fixed points. We conclude that $\widetilde F$ is trivial, which implies that $F$ itself is trivial. 
 \enskip 
 
 \medskip
Now, let us prove {\bf ii)}. We consider a closed connected orientable surface $M$ and two non trivial periodic orientation preserving homeomorphims $F_0$ and $F_1$ of period $q$ that are isotopic. We consider an isotopy $I=(F_t)_{t\in[0,1]}$ from $F_0$ to $F_1$. The sets ${\rm Fix}(F_0)$ and ${\rm Fix}(F_1)$ are finite by assertion {\bf i)}. Let us write $\widetilde M$ for the universal covering space of $M$ and $\cal T$ for the deck transformation group. Choose $z_0\in{\rm Fix}(F_0)$. Fix a lift $\widetilde z_0\in\widetilde M$ of $z_0$ and denote by $\widetilde F_0$ the lift of $F_0$ that fixes $\widetilde z_0$. One deduces that  $\widetilde F_0^{q}={\rm Id}_{\widetilde M}$. The isotopy $I$ may be lifted to an isotopy $\widetilde I=(\widetilde F_t)_{t\in[0,1]}$ starting from $\widetilde F_0$. It defines a lift   $\widetilde F_1$ of $F_1$. Of course, there exists $\tau_1\in {\cal T}$ such that $\widetilde F_1^q=\tau_1$. Write $\alpha_i: \tau\mapsto \widetilde F_i\circ \tau\circ\widetilde F_i^{-1}$ for the natural action of $\widetilde F_i$ on ${\cal T}$. Of course, one has $\alpha_0=\alpha_1$ because  $\widetilde F_0$ and $\widetilde F_1$ are the ends of the lifted isotopy $\widetilde I$.  We deduce that $\alpha_1^q=\alpha_0^q={\rm Id}_{\cal T}$, which means that $\tau_1$ belongs to the center of $\cal T$. But this center is trivial because the genus of $M$ is larger than one. We deduce that $\widetilde F_1^q={\rm Id}$. Consequently $\widetilde F_1$ has a unique fixed point $\widetilde z_1$. Replacing $\widetilde z_0$ with $\tau(\widetilde z_0)$ will change $\widetilde F_0$ into $\tau\circ \widetilde F_0\circ \tau^{-1}$ and $\widetilde F_1$ into $\tau\circ \widetilde F_1\circ \tau^{-1}$. So it will change $\widetilde z_1$ into $\tau(\widetilde z_1)$. Thus, we have defined a natural map from the fixed point set of $F_0$ to the fixed point set of $F_1$. Permuting $0$ and $1$, one gets a map from the fixed point set of $F_1$ to the fixed point set of $F_0$ which is reciprocal to the previous one. \enskip \hfill$\Box$

  \bigskip
  \bigskip

\noindent{\large {\bf 2. Acyclic and cellular decompositions of homeomorphisms of surfaces}}

\bigskip

\noindent{\bf General facts about upper semi-continuous decompositions}

\bigskip
An {\it upper semi-continuous decomposition} of a topological space $X$ is a partition $X=\bigsqcup_{i\in I} K_i$ in compact subsets such that every $K_{i_0}$ has a fundamental system of neighborhoods that  are {\it saturated} open sets (which means union of $K_i$). Equivalently, it means that the {\it saturation} of a closed set $Y$, which is the union of the $K_i$ that meet $Y$, is closed. In other terms, a partition $X=\bigsqcup_{i\in I} K_i$ in compact subsets is upper semi-continuous if the projection $\pi~: X\to I$ defined by $x\in K_{\pi(x)}$, is a closed map, when $I$ is munished with the quotient topology. We will write ${\cal D}=(K_i)_{i\in I}$ for such a decomposition, denote by ${\cal D}(x)$  the set $K_{\pi(x)}$ and more generally by ${\cal D}(Y)$ the saturation of a set $Y\subset X$. Recall some properties that will be useful later (see Whyburn \cite{W} for example)~:

\bigskip
\noindent{\sc Proposition 2.1~:}\enskip\enskip {\it Let  ${\cal D}=(K_i)_{i\in I}$ be an upper semi-continuous decomposition of a topological space $X$. 

\medskip
\noindent{\bf i)} If $X$ is Hausdorff, then $I$ is Hausdorff and $\pi$ is proper.

\medskip
\noindent{\bf ii)} If $X$ is Hausdorff and locally compact, then  $I$ is Hausdorff and locally compact.
 
 \medskip
 \noindent{\bf iii)} If $Y\subset X$ is saturated, the restriction of ${\cal D}$ to $Y$ defines an upper semi-continuous decomposition of $Y$. }

\bigskip
\noindent{\it Proof.} \enskip Let us prove {\bf i)}. Suppose that $i$ and $i'$ are two distinct elements of $I$. The space $X$ being Hausdorff, one can find two neighborhoods $U$ and $U'$ of $K_i$ and $K_{i'}$ respectively that are disjoint. The decomposition ${\cal D}$ being upper semi-continuous, one can find  two open saturated neighborhoods $V$ and $V'$ of $K_i$ and $K_{i'}$ respectively, such that $V\subset U$ and $V'\subset U'$. They project onto disjoint open sets of $I$ that contain $i$ and $i'$ respectively, because they are open, saturated and disjoint. We have proven that $I$ is Hausdorff. To prove that $\pi$ is proper we consider a compact set $I'\subset I$ and we want to prove that $\pi^{-1}(I')$ is compact. Let $(U_{\a})_{\a\in A }$ be an open cover of $\pi^{-1}(I')$. For every $i\in I'$, there exists a finite subcover $(U_{\a})_{\a\in A_i }$ of $K_i$. The decomposition being upper semi-continuous, one can find a saturated open set $V_i\subset \bigcup_{a\in A_i} U_a$ that contains $K_i$. We obtain an open cover $(\pi(V_i))_{i\in I'}$ of $I'$. By hypothesis, $I'$ is compact, so one can  find a finite subcover $(\pi(V_i))_{i\in I''}$ of $I'$. Observe now that $(U_{\a})_{\a\in A_i, i\in I''}$ is a finite subcover of $\pi^{-1}(I')$. We have proven that  $\pi^{-1}(I')$ is compact.

\medskip
To prove {\bf ii)} let us fix $i\in I$ and consider a neighborhood $O$ of $i\in I$. We want to  construct a compact neghborhood of $i$ which is contained in $O$. The set $\pi^{-1}(O)$ is a neighborhood of $K_i$ and $X$ is Hausdorff and locally compact. This implies that there exists an open cover $(V_{\a})_{\a\in A }$ of  $K_i$ by relatively compact saturated sets whose closure are all included in $\pi^{-1}(O)$. There exists a finite subcover $(V_{\a})_{\a\in A' }$ of $K_i$ because $K_i$ is compact. The compact set $\overline {\bigcup_{\a\in A'} V_{\a}}$ projects onto a compact neighborhood of $i$ that is included in $O$.

\medskip
The assertion {\bf iii)} is obviously true. \hfill$\Box$

\bigskip
\noindent{\sc Proposition 2.2~:}\enskip\enskip {\it Let $X$ and $Y$ be two Hausdorff topological spaces and suppose moreover than $Y$ is locally compact. If  $\pi~: X\to Y$ is a surjective, proper and continuous  map, then the decomposition $(\pi^{-1}(\{y\}))_{y\in Y}$ in fibers of $\pi$ is upper semi-continuous. }

\bigskip
\noindent{\it Proof.} \enskip Fix $a\in Y$ and denote by $(O_{\alpha})_{\alpha\in A}$ the family of compact neighborhoods of $a$. By hypothesis every set $\pi^{-1}(O_{\a})$ is compact. Moreover, one knows that $$\bigcap_{\a\in A} \pi^{-1}(O_{\a})= \pi^{-1}\left(\bigcap_{\a\in A} O_{\a}\right)=\pi^{-1}(\{a\}).$$ One deduces that for every neighborhood  $U$ of $ \pi^{-1}(\{a\})$, there exists a finite set $A'\subset A$ such that
$$\bigcap_{\a\in A'} \pi^{-1}(O_{\a}) \subset U.$$ Indeed one has
$$\bigcap_{\a\in A} \left(\pi^{-1}(O_{\a})\cap (X\setminus U)\right)=\emptyset.$$One gets a saturated open neighborhood of $\pi^{-1}(\{a\})$ that is contained in $U$ by considering $\pi^{-1}(O')$, where  $O'$ is an open neighborhood of $a$ that is incuded in $\bigcap_{\a\in A'} O_a$.\hfill$\Box$

 \bigskip 
 The  next result is a particular case of Proposition 2.2 and we will not be surprised to get a proof very similar.
  
\bigskip
\noindent{\sc Proposition 2.3~:}\enskip\enskip {\it Let $X$ be a locally compact Hausdorff space and $({\cal D}^j)_{j\in J}$ a family of upper semi-continuous decompositions on $X$, where ${\cal D}^j=(K_i^j)_{i\in I_j}$. Then the decomposition $\bigwedge_{j\in J}{\cal D}^{j}=(K_{\xi})_{\xi\in\Xi} $ is upper semi-continuous, where
$$\Xi=\left\{\xi=(i_j)_{j\in J} \in\prod_{j\in J} I_j,\vert \bigcap_{j\in J} K^j_{i_j}\not =\emptyset\right\}\enskip and\enskip
K_{\xi}=\bigcap_{j\in I} K^j_{i_j}.$$
We will call $\bigwedge_{j\in J}{\cal D}^{j}=(K_{\xi})_{\xi\in\Xi} $ the {\rm intersection} of the family $({\cal D}^j)_{j\in J}$.}

\bigskip
\noindent{\it Proof.} \enskip Observe first that each $K_{\xi}$ is compact. Fix $\xi=(i_j)_{j\in J}\in \Xi$ and a neighborhood $U$ of $K_{\xi}$. We want to construct an open neighborhood $W$ of $K_{\xi}$ that is included in $U$ and  $\bigwedge_{j\in J}{\cal D}^{j}$-saturated. By definition, one has $\bigcap_{ j\in J} K_{i_j}=K_{\xi}$, which implies that there exists a finite family $(j_k)_{1\leq k\leq r}$ such that $\bigcap_{ 1\leq k\leq r} K_{i_{j_k}}\subset U$. Consequently, one can write $ K_{i_{j_1}}\subset U\cup\left( X\setminus \left(\bigcap_{2\leq k\leq r} K_{i_{j_k}}\right)\right)$. The space $X$ being Hausdorff and locally compact, there exists a compact neighborhood $V_{1}$ of $K_{i_{j_1}}$ such that $ V_{1}\subset U\cup\left( X\setminus \left(\bigcap_{2\leq k\leq r} K_{i_{j_k}}\right)\right)$, which means that $V_{1}\cap\left(\bigcap_{2\leq k} K_{i_{j_k}}\right)\subset U$. A simple induction argument permits us to construct, for every $k\in\{1,\dots,r\}$, a compact neighborhood $V_{k}$ of $K_{i_{j_k}}$ such that  $\bigcap_{ 1\leq k\leq r} V_{k}\subset U$. Each decomposition ${\cal D}^{j_k}$, $1\leq k\leq r$, being upper semi-continuous, one can find, for every $k\in\{1,\dots,r\}$ a ${\cal D}^{j_k}$-saturated open neighborhood $W_k$ of $K_{i_{j_k}}$ that is included in $V_k$. The set $W=\bigcap_{1\leq k\leq r}W_k$ is a $\bigwedge_{j\in J}{\cal D}^{j}$-saturated open neighborhood of $K_{\xi}$ that is included in $U$.\hfill$\Box$

\bigskip
\noindent{\sc Proposition 2.4~:}\enskip\enskip {\it Let $X$ and $Y$ be two Hausdorff topological spaces, with $Y$ locally compact, and $\pi~: X\to Y$ a surjective, proper and continuous  map. If $F$ is a homeomorphim of $X$ that leaves $\pi$ invariant, which means that $\pi\circ F= \pi$, then every connected component of a fiber $\pi^{-1}(\{a\})$ has a fundamental system of neighborhoods that are open, connected,  union of connected components of fibers, and $F$-periodic.}

\bigskip
\noindent{\it Proof.} \enskip  Fix a connected component $K$ of $K_a=\pi^{-1}(\{a\})$ and a neighborhood $U$ of $K$. Denote by ${\cal V}$ the set of relatively compact open neighborhoods of $a$. For every $O\in{\cal V}$, write $W_O$ for the connected component of $\pi^{-1}(O)$ that contains $K$. The projection $\pi$ being proper, one knows that $W_O$ is a relatively compact open neighborhood of $K$. Each set $\overline{W_O}$ is connected and one has $W_{O\cap O'}\subset W_O\cap W_{O'}$. Consequently, one knows that $\bigcap_{O\in{\cal V}} \overline {W_{O}}$ is connected. This set, being included in $\pi^{-1}(\{a\})$, coincides with $K$. In particular, there exists $O\in{\cal V}$ such that $\overline{W_O}\subset U$.

\medskip

It remains to prove that $W_O$ is periodic. Observe that $\pi^{-1}(\{a\})$ is covered by finitely many connected components of $\pi^{-1}(O)$ and that $F$ induces a permutation on this set. \hfill$\Box$

\bigskip 
\noindent{\sc Corollary 2.5~:}\enskip\enskip {\it If $\cal D$ is an upper semi-continuous decomposition on a Hausdorff locally compact space $X$, then the decomposition in connected components of elements of $\cal D$ is also upper semi-continuous.}

\bigskip
\noindent{\it Proof.} \enskip  Write ${\cal D}=(K_i)_{i\in I}$. By Proposition 2.1, one knows that $I$ is Hausdorff, locally compact and that $\pi~: X\to I$ is surjective, proper and continuous. If one applies Proposition 2.4 to the identity map on $X$, one gets the corollary.\hfill$\Box$ 
 
\bigskip
 
We will say that an upper semi-continuous decomposition ${\cal D}=(K_i)_{i\in I}$ of a topological space is {\it monotone} if every $K_i$ is connected. The previous result asserts that every upper semi-continuous decomposition ${\cal D}=(K_i)_{i\in I}$ on a locally compact topological space induces naturally a monotone upper-semi continuous decomposition by taking the connected components of the $K_i$, we will called it the {\it induced monotone decomposition}.
If $ \pi': X\to I'$ is the projection, then one can write $\pi=\iota\circ \pi'$ where $\iota~: I'\to I$ is induced by the inclusion: $K'_i\subset K_{\iota(i')}$. The map $\iota$ is continuous because for every open set $O$ of $I$, the set $\pi^{-1}(O)=\pi'^{-1} (\iota^{-1}(O))$ is open in $M$, which means that $\iota^{-1}(O)$ is open in $I'$. The map $\iota$ is proper because, for every compact set $K\subset I$, one can write $\iota^{-1}(K)=\pi'(\pi^{-1}(K))$, which implies that $\iota^{-1}(K)$ is compact.

\medskip

Observe that if ${\cal D}=(K_i)_{i\in I}$ is a monotone upper semi-continuous decomposition on a topological space, then the inverse image of every connected set $I'\subset I$ by the projection $\pi~: X\to I$ is connected. Indeed, if $\theta~:\pi^{-1}(I')\to\{0,1\}$ is a continuous map, it is constant on every set $K_i$, $i\in I'$, because such a set is connected. So, one can write $\theta=\Theta\circ\pi$, where $\Theta~: I'\to \{0,1\}$ is continuous. The fact that $I$ is connected implies that $\Theta$ is constant. In particular, the projection $\pi~: X\to I$ induces a bijection between the set of connected components of  
$X$ and the set of connected components of $I$

\bigskip
\bigskip

\noindent{\bf Upper semi-continuous decompositions and dynamics}

\bigskip
Let $X$ be a Hausdorff locally compact space, $F$ a homeomorphism of $X$ and ${\cal D}=(K_i)_{i\in I}$ an upper semi-continuous decomposition of $X$. Let us say that $F$ {\it acts on } ${\cal D}$ if there is a bijection $\varphi~: I\to I$ such that $F(K_i)=K_{\varphi (i)}$ for every $i\in I$. In that case,  $\varphi$ is a homeomorphism of the space $I$ munished with the quotient topology,  and the projection $\pi~:X\to I$ is a factor map from $F$ to $\varphi$. In the case where $\varphi$ is the identity map, that means  if every $K_i$ is invariant by $F$, let us say that ${\cal D}$ is {\it invariant} (or $F$-{\it invariant}); in the case where ${\cal D}$ is invariant by a power $F^q$ of $F$, which means that $\varphi^q={\rm Id}_I$, let us say that ${\cal D}$ is {\it periodic} (or {\it $F$-periodic}). If $F$ acts on  ${\cal D}=(K_i)_{i\in I}$, then $F$ acts naturally on the induced monotone  decomposition ${\cal D}'=(K_i)_{i\in I}$
and the induced map $\varphi'$ is an extension of $\varphi$ with $\iota$ as a factor map. Noting that the image by $\pi'$ of a connected set of $X$ is a connected set of $I'$ by continuity of $\pi'$, observe that Proposition 2.4 implies immediatly the following result (and is equivalent to it in the case where $X$ is locally compact):

\bigskip
\noindent{\sc Proposition 2.6~:}\enskip\enskip {\it Let $X$ be a  Hausdorff locally compact space,  $F$ a homeomorphism of $X$ and ${\cal D}$ an upper semi-continuous decomposition of $X$ that is invariant by $F$. Let denote by ${\cal D}'=(K_i')_{i'\in I'}$ the induced monotone decomposition and by $\varphi'$ the induced map on $I$. Then every point $i'\in I'$ has a fundamental system of neighborhoods that are open, connected and $F$-periodic.}

\bigskip
\bigskip

\noindent{\bf Local decompositions}

\medskip
Let $X$ be a Hausdorff locally compact space and $\overline X$ the Alex\-an\-drov compactification of $X$. A {\it local upper semi-continuous decomposition}  is a couple  ${\cal D}=\left(U,(K_i)_{i\in I}\right)$, where $U$ is an open subset $U$ of $X$, called the {\it domain} of $\cal D$, and $U=\bigsqcup_{i\in I} K_i$ an upper semi-continuous decomposition on $U$. By Proposition 2.1, one knows that $I$ is Hausdorff and that the projection $\pi~:U\to I$ is continuous and proper. One deduces that the saturation ${\cal D}(K)=\pi^{-1}(\pi(K))$ of any compact set $K$ of $U$ is compact. Consequently, one gets an upper semi-continuous decomposition $\overline{\cal D}$ of $\overline X$, the {\it augmented decomposition}, by adding the compact set $\overline X\setminus U$. 

\medskip
Let us introduce another useful notion. We will say that  ${\cal D}=\left(U',(K_{i'})_{i'\in I'}\right)$ is {\it  finer} than ${\cal D}=\left(U,(K_i)_{i\in I}\right)$ if $U\subset U'$ and if every $K_i$, $i\in I$, is ${\cal D}'$-saturated.  If $(D^j)_{j\in J}$ is a family of local upper semi-continuous decompositions, where ${\cal D}^j=\left(U^j,(K^j_i)_{i\in I}\right)$, we will denote by $\bigwedge_{j\in J}{\cal D}^j$ the decomposition whose domain is  $\bigcup_{j\in J} U^j$, defined as the restriction to  $\bigcup_{j\in J} U^j$ of $\bigwedge_{j\in J}\overline{{\cal D}^j}$. It is finer than every decomposition  ${\cal D}^j$, $j\in J$. Here again, we call it the intersection of the ${\cal D}^j$.

\medskip
If $F$ is a homeomorphim of $X$, let us say that $F$ acts on a local upper semi-continuous decomposition ${\cal D}=\left(U,(K_i)_{i\in I}\right)$ if $U$ is invariant by $F$ and if $F\vert_U$ acts on ${\cal D}$. If the induced map $\varphi$ on $I$ is the identity, we will say that ${\cal D}$ is invariant; if it periodic, we will say that ${\cal D}$ is periodic.

\bigskip
\bigskip

\noindent{\bf Acyclic decompositions}

\medskip

We will focus now on decompositions on surfaces. Until the end of the section, we mean by $M$ an orientable surface without boundary, not necessarily compact, not necessarily connected. Let us introduce some definitions.

\bigskip
\noindent{\sc Proposition 2.7~:}\enskip\enskip {\it Let $U$ be a connected open set of a surface $M$. The following conditions are equivalent:

\smallskip
\noindent{\bf i)}\enskip\enskip $U$ is included in an open disk of $M$ or the connected component of $M$ that contains $U$ is the $2$-sphere;

\smallskip
\noindent{\bf ii)}\enskip\enskip every loop in $U$ is homotopic to zero in $M$;

\smallskip
\noindent{\bf iii)}\enskip\enskip every simple closed curve in $U$  borders a disk in $M$.

\smallskip
\noindent If these properties are not satified we will say that $U$ is {\rm essential } }

\bigskip
\noindent{\it Proof.} \enskip   The fact that {\bf i)} implies {\bf ii)} and {\bf iii)} is obvious. Let us prove that {\bf ii)} implies {\bf i)}. Suppose that the connected component $M'$ of $M$ that contains $U$ is not the $2$-sphere. Write $\widetilde M'$  for its universal covering space and $\Pi~:\widetilde M'\to M'$ for the covering projection. Fix a connected component $\widetilde U$ of $\Pi^{-1}(U)$. The assertion {\bf ii)} tells us that $\widetilde U$ does not meet its images by the deck transformations and that $\Pi$ induces a homeomorphism between $\widetilde U$ and $U$. Let $\widetilde V$ be the union of $\widetilde U$ and of all the compact connected components of $\widetilde M'\setminus \widetilde U$. One gets an open disk that does not meet its images by the deck transformations and that projects into an open disk $V$ that contains $U$. Let us prove now that {\bf iii)} implies {\bf ii)}. Let $\Gamma$ be a loop in $U$. One can find a connected compact surface with boundary $N\subset U$ that contains $\Gamma$. Every boundary circle of $N$ borders a disk in $M$. If this disk  contains $\Gamma$  we are done. So one can suppose that every boundary circle of $N$ borders a disk disjoint from the interior of $N$. Consequently $M'$ is obtained from $N$ by pasting a disk on every boundary circle of $N$. If {\bf iii)} is satisfied, the only possibility is that $M'$ is a $2$-sphere, which implies {\bf ii)}. \hfill$\Box$

\bigskip
We will say that an open set $U$ of $M$ is {\it essential} if every connected component of $U$ is essential. We will say that a compact set $K\subset M$ is {\it cellular} if there exists a non increasing sequence of closed disks $(D_n)_{n\geq 0}$ such that $K=\bigcap_{n\geq 0} D_n$. We will say that a compact set $K\subset M$ is {\it acyclic} is every connected component of $K$ is cellular. A classical result of plane topology tells us that if $U\subset M$ is an open disk, then a compact set $K\subset U$ is acyclic if and only if $U\setminus K$ is connected As a consequence, in the case where $M$ is connected, a compact set $K$ which is contained in an open disk $U$ and such that $M\setminus K$ is connected is  acyclic. Indeed, there is no connected component of $U\setminus K= U\cap(M\setminus K)$ that is relatively compact in $K$ (it is an immediate consequence of a classical theorem of topology about connectedness, see for example Hocking-Young \cite{HoYou}, Theorem 2-16).

\medskip

We will call {\it essential decomposition} every local upper semi-continuous decomposition whose domain is essential. We will call {\it acyclic decomposition} every local upper semi-continuous decomposition  ${\cal D}=\left(U,(K_i)_{i\in I}\right)$ such that all the $K_i$ are acyclic. If all the $K_i$ are cellular, we will say that ${\cal D}$ is  a {\it cellular decomposition}. In other words a cellular decomposition is a monotone acyclic decomposition. In particular, the monotone decomposition induced by an acyclic decomposition is a cellular decomposition. 

Let us state one of the fundamental result of this section:

\bigskip
\noindent{\sc Proposition 2.8~:}\enskip\enskip {\it If $({\cal D}^j)_{j\in J}$ is a family of essential acyclic decompositions of $M$, then $\bigwedge_{j\in J}{\cal D}^j$ is an essential acyclic decomposition.}

\bigskip
\noindent{\it Proof.} \enskip The fact that $\bigwedge_{j\in J}{\cal D}^j$ is essential is obvious because its domain is the union of the domains of the ${\cal D}^j$. Let us prove that it is acyclic. Write ${\cal D}^j=\left(U^j,(K^j_i)_{i\in I^j}\right)$ and $\overline{{\cal D}^{j}}=\left(\overline M,(K^j_i)_{i\in {\overline {I^j}}}\right)$, where $\overline {I^j}=I^j\sqcup\{\infty\}$ and $K^j_{\infty}=\overline {M}\setminus U^j$.  We must prove that every element $K=\bigcap_{j\in J} K^j_{i_j}$ of the decomposition $\bigwedge_{j\in J}\overline{{\cal D}^j}$, which is not $\bigcap_{j\in J} K^j_{\infty}$, is acyclic. That means that every connected component of $K$ is cellular. Fix such a component $X$ and write $M'$ for the connected component of $M$ that contains $X$. If $i_j\not=\infty$ write $X^j_{i_j}$ for the connected component of $K^j_{i_j}$ that contains $X$; if $i_j=\infty$ write $X^j_{i_j}=M'\setminus U_j$.  Then consider the closed set $X'=\bigcap_{j\in J} X^j_{i_j}$. From the inclusions $X\subset X'\subset K$, one deduces that $X$ is a connected component of $X'$ and so it is sufficient to prove that $X'$ is acyclic. Observe that  $X^j_{i_j}$ is cellular if $i_j\not=\emptyset$ because $K^j_{i_j}$ is acyclic. Moreover, there exists  $j_0\in J$ such that $i_{j_0}\not =\infty$. One deduces that there exists an open disk $U\subset M'$ such that $X'\subset X^{j_0}_{i_{j_0}}\subset U$. Therefore, to prove that $X'$ is acyclic, it is sufficient to prove that $M'\setminus X'$ is connected.  By hypothesis, one knows that $M'\not=S^2$ because there is no essential open set on $S^2$. Let us write $$\eqalign{M'\setminus X'&=\left(\bigcup_{j\in J\,\vert\,i_j\not=\infty} M'\setminus X^j_{i_j}\right)\cup\left(\bigcup_{j\in J\,\vert\, i_j=\infty}  M'\setminus X^j_{i_j}\right),\cr
&=\left(\bigcup_{j\in J\,\vert\,i_j\not=\infty} M'\setminus X^j_{i_j}\right)\cup\left(\bigcup_{j\in J\,\vert\, i_j=\infty} M'\cap U^j\right).\cr}$$

In the previous formula, if $i_j\not=\infty$, one has $X^j_{i_j}\subset M'$ and so  $M'\setminus X^j_{i_j}$ is connected because $X^j_{i_j}$ is cellular. A surface that is not a sphere cannot be covered by two open disks, and so cannot be covered by two cellular compact sets. One deduces that $M'\not=X^{j_0}_{i_{j_0}}\cup X^{j}_{i_j}$ if $i_j\not=\infty$, which means that the connected set  $M'\setminus X^j_{i_j}$ intersects the connected set $M'\setminus X^{j_0}_{i_{j_0}}$. Every set $U^j$, $i_j=\infty$, being essential, none of its connected components is included in $X^{j_0}_{i_{j_0}}$. In other words, every connected component of $U^j$ that is contained in $M'$ meets $M'\setminus X^{j_0}_{i_{j_0}}$. We have written $M'\setminus X'$ as the union of connected open sets which all meet $M'\setminus X^{j_0}_{i_{j_0}}$. This last set being connected and included in $M'\setminus X'$, this implies that $M'\setminus X'$ itself is connected.\hfill$\Box$

\bigskip
\bigskip
\noindent{\bf Induced cellular decompositions}

\bigskip
 Let $F$ be a homeomorphism of $M$. If $(D^j)_{j\in J}$ is a family of invariant local upper semi-continuous decompositions, it is clear that $\bigwedge_{j\in J}{\cal D}^j$ is also invariant. By Proposition 2.8, one can say that the intersection of all the invariant essential acyclic decompositions is itself acyclic, it is the finest among all invariant essential acyclic decompositions, we call it the {\it canonical invariant essential acyclic decomposition} of $F$. 

\medskip

Let us state now the  main result of Section 2. In particular, it asserts that the monotone decomposition induced by the canonical invariant essential acyclic decomposition of $F$ is periodic on each connected component of its domain.

\bigskip
\noindent{\sc Proposition 2.9}\enskip\enskip {\it Let ${\cal D}=(K_i)_{i\in I}$ be a global acyclic decomposition on a connected surface $M$ that is invariant by a homeomorphism $F$. Then the cellular decomposition naturally induced by ${\cal D}$ is $F$-periodic.}

\bigskip
\noindent{\it Proof.} \enskip  Write ${\cal D}'=(K_{i'})_{i'\in I'}$ for the induced cellular decomposition and $\varphi'$ for the bijection of $I'$ induced by $f$. According to a result of Moore \cite{Moo}, one knows that the topological space $I'$ is a surface homeomorphic to $M$. Proposition 2.6 tells us that every point $i'\in I'$ has a fundamental system of neighborhoods that are open, connected and $\varphi'$-periodic.  According to Theorem 1.1, one concludes that $\varphi'$ is periodic. \hfill$\Box$

\bigskip 
\noindent{\bf Remark:}\enskip In fact, Moore's Theorem is not necessary. Indeed, what has been proved above is the fact that every connected component of a $K_i$ is cellular and has a fundamental system of neighborhoods that are open, connected, union of such components, and $F$-periodic. The proof of  Theorem 1.1 may be slightly changed and adapted to this situation to prove that there is a power of $F$ that fixes all these connected components.

\bigskip
Let us continue with the properties of acyclic and cellular decompositions.

\bigskip
\noindent{\sc Theorem 2.10}\enskip\enskip {\it Let $F$ be a homeomorphism of a surface $M$.  The intersection of all the invariant essential cellular decompositions is itself an invariant essential cellular decomposition, it is the finest among all such decompositions, we call it the {\rm canonical invariant essential cellular decomposition} of $F$.}

\bigskip
\noindent{\it Proof.} \enskip Of course, one can assume that there exists at least one invariant essential cellular decomposition, otherwise the result is true with  an invariant essential cellular decomposition of empty domain. According to Proposition 2.8, the intersection of all the invariant essential cellular decompositions is an invariant essential acyclic decomposition ${\cal D}=(U, (K_{i})_{i\in I})$. We consider the induced cellular decomposition ${\cal D}'=(U, (K'_{i'})_{i'\in I'})$. To prove Proposition 2.8, it is sufficient to prove that ${\cal D}'$ is invariant.  Indeed ${\cal D}'$  will coincide with ${\cal D}$ because each of these two decompositions will be finer than the other one. In other words, one must prove that the homeomorphism $\varphi'$ that is naturally induced by $F$ on $I'$ is equal to the identity.  Let $I'_*$ be a connected component of $I'$. Recall that the preimage of $I'_*$ by the projection $\pi': U\to I'$ is a connected component $U_*$ of $U$ because ${\cal D}'$ is monotone. Fix $z_*\in U_*$.  By definition of $U$, there exists an invariant essential cellular decomposition  ${\cal D}''=(U'', (K''_{i''})_{i\in I''})$ whose domain contains $z_*$. Write $\pi'' : U''\to I''$ for the projection and  consider a connected compact neighborhood $X''_*$ of $i''_*=\pi''(z_*)$ in the surface $I''$. The set $\pi''^{-1}(X''_*)$ is compact because $\pi''$ is proper and connected because ${\cal D}''$ is monotone. Therefore $\pi''^{-1}(X''_*)$ is included in $U_*$. One deduces that $X'_*=\pi'(\pi''^{-1}(X''_*))$ is a compact and connected subset of $I'_*$. For every $i''\in X''_*$, there exists a fixed point $z''$ of $F$ in $K''_{i''}$ by Cartwright-Littlewood Fixed Point Theorem and $\pi'(z'')$ is a fixed point of $\varphi'$. It is clear that points of $U''$ that project on different points of $I''$ also project on different points of $I'$. So we have proved that the compact set $X'_*$ contains infinitely many fixed points of $\varphi'$ which implies, by Proposition 1.6, that $\varphi'$ fixes all the points of $I''_*$. \hfill$\Box$

 \bigskip
\noindent{\bf Remark:}\enskip Observe that the domain $U$ of the invariant essential acyclic decomposition ${\cal D}$ of $F$ can be written $U=\cup_{k\geq 1}U^k$, where $U^k$ is the domain of the invariant essential cellular decomposition of $F^k$. Indeed, the map $\varphi'$ naturally defined by $F$ on the cellular decomposition induced by ${\cal D}$ satifies the property {\bf (P)} of Theorem 1.1. Thus, every connected component $V$ of $U$ is fixed by a power $F^k$, $k\geq 1$, and the map $\varphi'^k$ defined on the restricted cellular decomposition is periodic by Proposition 2.9. Consequently, there exists $k'\geq 1$ such that $V\subset U^{kk'}$. Conversely, for every $k\geq 1$, the map $F$ acts on the invariant essential cellular decomposition $(U^k, {\cal D}^k)$ of $F^k$, because $F$ commutes with $F^k$. Write ${\cal D}^k=(K_i)_{i\in I_k}$. The domain $U^k$ is invariant by $F$ and every set $\cup_{0\leq j<k} F^j(K_i)$ is invariant by $F$ 	and acyclic, being the union of $d$ disjoint cellular sets, where $d$ divides $ k$. It remains to observe that the family of such sets defines a semi-continuous decomposition of $U^k$. Indeed for every open neighborhood $W$ of $\cup_{0\leq j<k} F^j(K_i)$, one can find a  ${\cal D}^k$-saturated neighborhood $O$ of $K_i$ such that $\cup_{0\leq j<k} F^j(O)\subset W$.

\bigskip
\noindent{\sc Proposition 2.11.}\enskip\enskip {\it  Let $F$ be a homeomorphism of a surface $M$. Write ${\cal D}=(U, (K_i)_{i\in I})$ for its canonical invariant essential acyclic decomposition and ${\cal D}'=(U, (K'_{i'})_{i'\in I'})$ for the induced cellular decomposition. Denote by $\pi':U\to I'$ the projection and by $\varphi'$ the map induced on $I'$ by $F$. If $U'$ is the union of the components $W$ of $U$ that are invariant by $F$ and such that $\varphi'\vert_{\pi'(W)}$ is the identity, then the canonical invariant cellular decomposition of $F$ is the restriction of $\cal D$ to $U'$.}

\bigskip
\noindent{\it Proof.} \enskip  The decomposition ${\cal D}$ is obviously finer than the canonical invariant cellular decomposition ${\cal D}''=(U'', (K''_{i''})_{i''\in I''})$ of $F$. This implies that $U''\subset U$. To prove the proposition, it is sufficient to prove that, for every component $W$ of $U$ that meets $U''$, one has $\varphi'_{W}={\rm Id}_{\pi'(W)}$. Indeed, we will deduce that $W\subset U''$. The argument is the same as in the previous result: the fact that $W$ meets $U''$ implies first that $W$ is invariant by $\varphi'$, then that the set of fixed point of $\varphi'_W$ is not discrete. The map $\varphi'_W$ being periodic one deduces that $\varphi_W={\rm Id}_{\pi'(W)}$.\hfill$\Box$

\bigskip
\noindent{\sc Corollary 2.12}\enskip\enskip {\it Let $F$ be a homeomorphism of a surface $M$ and $U$ the domain of the canonical invariant cellular (or acyclic) decomposition ${\cal D}=(U, (K_i)_{i\in I})$ of $F$. Then every simple closed curve in $U$ that is homotopic to zero in $M$ borders a closed disk $D\subset U$. }

\bigskip
\noindent{\it Proof.} \enskip  According to the remark that follows Therem 2.10, it is sufficient to prove the result for the celullar decomposition. Of course, one can suppose that $M$ is connected and that $U\not=\emptyset$, which implies that $M\not=S^2$. Every simple closed curve $C\subset U$ homotopic to zero in $M$ borders a unique closed disk $D\subset M$. We want to prove that this disk is included in $U$. The saturation ${\cal D}(C)$ is an invariant compact connected set and $F$ induces a natural bijection on the set of connected components of $M\setminus {\cal D}(C)$. This bijection is the identity. Indeed, every connected component $W$ of $M\setminus {\cal D}(C)$ meets $U$ because it is open and its frontier is included in ${\cal D}(C)$. This implies that $W$ meets a set $K_i$. But $K_i$ being connected and disjoint from ${\cal D}(C)$ should be contained in $W$. As we know that $K_i$ contains a fixed point of $F$, by Cartwright-Littlewod Fixed Point Theorem, we deduce that $W$ is $F$-invariant. To prove that $D$ is included in $U$, one must prove that every connected component of $D\setminus {\cal D}(C)$ is included in $U$. Observe that such a component is an open disk because its complement in the interior of $D$, has no compact connected component.

Suppose that it is not the case, that means suppose that there exists a connected component $V$ of $D\setminus {\cal D}(C)$ that is not included in $U$. Let us explain why  $V\setminus U$ is an acyclic compact set invariant by $F$. The fact that $V\setminus U$ is compact comes from the inclusions ${\rm Fr }(V)\subset{\cal D}(C)\subset U$. The set $V\setminus U$ is invariant because  $V$ and $U$ are invariant. It remains to prove that $V\setminus U$ is acyclic. Of course it is included in a disk, the disk $V$ itself. So, to conclude, one must prove that  the complement of $V\setminus U$ in $M$ is connected. This arises from the connectedness of $M\setminus D$ and of the fact that there is no connected component of $U$ that is contained in $V$. Indeed, $U$ is supposed to be essential.

To find a contradiction, let us add $V\setminus U$ to $\cal D$ to get an invariant acyclic decomposition finer than ${\cal D}$. One of the connected component of its domain, the one that contains $V$, meets the domain of $\cal D$, without being included in its domain. This contredicts Proposition 2.9. \hfill$\Box$

 \bigskip
\noindent{\bf Remark 1:}\enskip In the case where $M$ is a closed surface, or a surface of finite type, one deduces that the domain $U$ of the  canonical invariant cellular (or acyclic) decomposition has finitely many ends (if it is not empty) and that these ends are trivial in the following sense: one can write $U=\bigcup_{n\geq 0} Q_n$ where $(Q_n)_{n\geq 0}$ is a sequence of  compact surfaces with boundaries, such that
 $Q_{n}\subset{\rm Int }(Q_{n+1})$, and such that $Q_{n+1}\setminus{\rm Int}(Q_n)$ is a union of essential compact annuli.
 
 \bigskip

\noindent {\bf Remark 2:}\enskip Let $M$ be a surface and $\cal F$ a subgroup of the group of homeomorphisms of $M$. One can define the canonical invariant essential acyclic decomposition of $\cal F$ as the intersection of all the essential acyclic decompositions that are invariant by every element of $\cal F$. Similarly, one can define the canonical invariant essential cellular decomposition of $\cal F$. It is the finest among all the invariant essential cellular decompositions that are invariant by every element of $\cal F$. Observe that its domain also satisfies the conclusion of Corollary 2.12.

 \bigskip
\noindent{\bf Example.} \enskip\enskip Consider the homeomorphism
$$\eqalign{ F~: \R/\Z\times\R&\to\R/\Z\times\R\cr (\theta,r)&\mapsto(\theta+\psi(r), r)\cr}$$where
 $\psi~:\R\to\R$ is continuous. For every $s\in\Q$, set $J_s={\rm Int}(\psi^{-1}(\{s\}))$. Define $$J={\rm Int}(\varphi^{-1}(\Z))=\bigcup_{s\in \Z} J_s$$ and $$J'={\rm Int}(\psi^{-1}(\Q))=\bigcup_{s\in \Q} J_s.$$ One can verify that the domain of the canonical invariant essential cellular decomposition is $U=\R/\Z\times J={\rm Int}({\rm Fix}(F))$ and that the decomposition on this set is the partition in points. Similarly, the domain of the canonical invariant essential acyclic decomposition is $U'=\R/\Z\times J'={\rm Int}({\rm Per}(F))$ and the decomposition is the decomposition in orbits of $F$. On a given connected component of $U'$, all the orbits have the same cardinal.

 \medskip
 This example permits to understand why looking at essential decompositions is fundamental in this section. If one adds the upper end $\infty$ to this annulus one gets a homeomorphism of a plane with no finest invariant cellular upper semi-continuous decomposition. Indeed, the previous cellular decomposition defines a cellular upper semi-continuous decomposition of this new map, but there is a lot of other such decompositions and none of them are comparable. To be more precise suppose that $\psi$ is increasing and that the $J_s$, $s\in\Z$, have non empty interior. For every $k\in\Z$ we have a decomposition whose domain is $\left({\bf R}/{\bf Z}\times \left(I\cup [\max J_k, +\infty)\right)\right)\cup\{\infty\}$ that consists of the disk $\left({\bf R}/{\bf Z}\times[\max J_k, +\infty)\right)\cup\{\infty\}$ and of the points of ${\bf R}/{\bf Z}\times \bigcup_{s\leq k} J_s$. In fact, this example is enlightening to understand why the domains must be essential in the statement of Proposition 2.8. Indeed add the lower end $-\infty$ to get a sphere. Fix $r$ in the interior of $J_0$. We have a decomposition whose domain is $\left({\bf R}/{\bf Z}\times \left(J\cup [r, +\infty)\right)\right)\cup\{\infty\}$ that consists of the disk $\left({\bf R}/{\bf Z}\times[r, +\infty)\right)\cup\{\infty\}$ and of the points of ${\bf R}/{\bf Z}\times \left(J\cap(-\infty, r)\right)$.  We have another decomposition whose domain is $\left({\bf R}/{\bf Z}\times \left(J\cup (-\infty, r]\right)\right)\cup\{-\infty\}$ that consists of the disk $\left({\bf R}/{\bf Z}\times(-\infty, r]\right)\cup\{-\infty\}$ and of the points of ${\bf R}/{\bf Z}\times \left(J\cap(r,+\infty)\right)$.  Both of them are invariant cellular decompositions but the intersection is not acyclic, because one of the elements of this decomposition is the circle  ${\bf R}/{\bf Z}\times \{r\}$.

\bigskip
\bigskip

\noindent{\large {\bf 3. Local Anosov homeomorphisms and centralizers}}

\bigskip
\noindent{\bf Hyperbolic automorphisms of $\T^2$ and centralizers}

\bigskip

We begin this subsection by recalling some well known facts about hyperbolic linear automorphisms of the torus based on the Shadowing Lemma, which are due to Franks \cite{Fr}, and will state a much recent one that is proved in \cite{Mark} (see also \cite{FrHa2}). Then it will be easy to copy the various proofs in the more general situation of local Anosov homeomorphisms. Let us fix a hyperbolic linear automorphism $L$ of $\T^2$, which preserves the orientation and has a unique fixed point, for example consider
$$L: (x,y)\mapsto (2x+y, x+y).$$

Denote by 
$l~: (x,y)\mapsto (2x+y, x+y)$
 the linear automorphism of $\R^2$ that lifts $L$ and fix a norm $\Vert \enskip\Vert$ on $\R^2$. Say that a sequence $(z_k)_{k\in\Z}$ in $\R^2$ is a {\it $C$-pseudo-orbit} of $l$ if for every $k\in\Z$, one has $\Vert l(z_k)-z_{k+1}\Vert \leq C$. Let us recall the {\it Shadowing Lemma}:

\bigskip
\noindent{\sc Theorem 3.1~:}\enskip\enskip {\it There exists a constant $C_*>0$ satisfying the following property:  for every $C$-pseudo-orbit $(z_k)_{k\in\Z}$, there exists a unique point $z'\in \R^2$ such that the sequence $\left(z_k-l^k(z')\right)_{k\in\Z}$ is bounded and moreover one has $\Vert z_k-l^k(z')\Vert\leq C_*C$.}

 \bigskip
 Let us state now some consequences of the Shadowing Lemma.
 
 \bigskip
\noindent{\sc Proposition 3.2~:}\enskip\enskip {\it The unique continuous transformation of $\T^2$ that is homotopic to the identity and commutes with $L$ is the identity.}

\bigskip
\noindent{\it Proof.} \enskip  Suppose that $H$ commutes with $L$ and is homotopic to the identity. The fact that $L$ has a unique fixed point implies that there exists a unique lift of $H$ to $\R^2$ that commutes with $l$. Indeed, fix a lift $h$ of $H$ to $\R^2$ and observe that there exists $\omega\in\Z^2$ such that $l\circ h=t_{\omega}\circ h\circ l$, where $t_{\omega}~: x\mapsto x+{\omega}$ is the translation of vector ${\omega}\in\Z^2$. To find a lift of $H$ that commutes with $l$, one must find ${\omega}'\in\Z^2$ such that $t_{{\omega}'}\circ h\circ l=l\circ t_{{\omega}'}\circ h$. As we have
$$l\circ t_{{\omega}'}\circ h=t_{l({\omega}')}\circ l\circ h=t_{l({\omega}')+{\omega}}\circ h\circ l,$$ we must solve the equation $${\omega}'=l({\omega}')+{\omega}.$$The fact that $L$ has a unique fixed point means that the degree of $L-{\rm Id}_{\T^2} $ is equal to $\pm1$,  or equivalently that ${\rm det}(l-{\rm Id}_{\R^2})=\pm 1$. Thus, the previous equation has a unique integer solution. Let us suppose now that $h$ commutes with $l$. The map $h-{\rm Id}_{\R^2}$ is bounded because it is $\Z^2$-periodic. One deduces that the sequence $$\left( l^k(h(z))- l^k(z)\right)_{k\in\Z}= \left( h(l^k(z))- l^k(z)\right)_{k\in\Z}$$ is bounded: the $l$-orbit  $h(z)$ shadows the $l$-orbit of $z$. By uniqueness of the shadowing orbit given by Theorem 3.1, one deduces that $h(z)=z$, which implies that $H$ is the identity. \hfill$\Box$

\bigskip
\noindent{\sc Proposition 3.3~:}\enskip\enskip {\it Every homeomorphism $G$ of $\T^2$ isotopic to $L$ is a continuous extension of $L$. More precisely, there exists a unique factor map $H~: \T¬^2\to\T^2$ from $G$ to $L$ that is homotopic to the identity.}

\bigskip \noindent{\it Proof.}\enskip Let $(G_s)_{s\in[0,1]}$ be an isotopy in $\T^2$ from $L$ to $G$. One can lift this isotopy on $\R^2$ to an isotopy $(g_s)_{s\in[0,1]}$ such that $g_0=l$. Of course, one has $g_s\circ t_{\omega}=t_{l({\omega})}\circ g_s$, for every ${\omega}\in \Z^2$. Observe that there exists $C>0$ such that for every $s\in[0,1] $ and every $z\in \R^2$, one has
$\Vert g_s(z)-l(z)\Vert \leq C$. This implies that the sequence $( g_s^k(z))_{k\in\Z}$ is a $C$-pseudo-orbit of $l$~: for every $k\in\Z$, one has 
$$ \Vert g_s^{k+1}(z) - l(g_s^{k}(z))\Vert \leq C.$$ By Theorem 3.1, there exists a unique point $z'\in \R^2$ such that the sequence $\left(g_s^k(z)-l^k(z')\right)_{k\in\Z}$ is bounded. Moreover, this sequence is bounded by $C_*C$. We denote by $h_s(z)$  this point. The graph of the function $(s,z)\mapsto h_s(z)$ is the following closed set
$$\left\{(s,z,z')\in[0,1]\times\R^2\times \R^2\enskip\vert\enskip k\in\Z\Rightarrow\Vert  g_s^k(z)-l^k(z')\Vert \leq C_*C\right\}.$$
Applying the Closed Graph Theorem on every compact set of $[0,1]\times \R^2$, one can affirm that the map $(s,z)\mapsto h_t(z)$ is continuous on $[0,1]\times\R^2$. By uniqueness, one gets

\smallskip
\noindent-\enskip\enskip  $h_0={\rm Id}_{\R^2}$;

\smallskip
\noindent-\enskip\enskip   $h_s\circ g_s=l\circ h_s$;

\smallskip
\noindent-\enskip\enskip   $h_s\circ t_{{\omega}}=t_{{\omega}}\circ h_s$ for every  $\omega\in\Z^2$.

\smallskip
Therefore, $(h_s)_{s\in[0,1]}$ lifts a continuous family $(H_s)_{s\in[0,1]}$ of transformations of $\T^2$, such that $H_0={\rm Id}$ and such that $H_s\circ G_s=L\circ H_s$ for every $s\in[0,1]$. The map $H_0$ inducing the identity on the group of homology $H_2(\T^2,\Z)$, it is the same for $H=H_1$, which implies that $H$ is onto.

 \medskip
To conclude the proof, it remains to prove the unicity of the factor map $H$. One must prove that every continuous transformation $H'$ of $\T^2$ homotopic to the identity such that $H'\circ G=L\circ H'$ coincides with $H$. Consider a homotopy $(H'_s)_{s\in[0,1]}$ from ${\rm Id}_{\T^2}$ to $H$. Lift it to a homotopy $(h'_s)_{s\in[0,1]}$ such that $h'_0={\rm Id}_{\R^2}$ and set $h'=h'_1$. Observe that  $h'\circ t_{\omega}=t_{\omega}\circ h'$ for every ${\omega}\in\Z^2$.  As a consequence of the equality $L\circ H'=H'\circ G$, one knows that there exists ${\omega}_0\in\Z^2$ such that $l\circ  h'=t_{{\omega}_0}\circ h'\circ g$. By hypothesis, $L$ has a unique fixed point and ${\rm det}(l-{\rm Id})=\pm 1$. One deduces that there exists ${\omega}_1\in\Z^2$ unique, such that $l({\omega}_1)-{\omega}_1=-{\omega}_0$. Set $h''=t_{{\omega}_1}\circ h'$ and observe that
$$l\circ h''=l\circ t_{{\omega}_1}\circ h'=t_{l({\omega}_1)}\circ l\circ h'=t_{l({\omega}_1)}\circ t_{{\omega}_0}\circ h'\circ g=t_{{\omega}_1}\circ h'\circ g=h''\circ g.$$ One deduces that $l^k\circ h''=h''\circ g^k$ for every $k\in\Z$. The map $h''-{\rm Id}_{\R^2}$ being $\Z^2$-periodic is bounded. This implies that $h''(z)$ is the unique point $z'$ such that the sequence $\left( g^k(z)-l^k(z')\right)_{k\in\Z}$ is bounded, that means $h''(z)=h(z)$. One deduces that $H'=H$. \hfill$\Box$
\bigskip

The next result is not so familiar (and is implicitely proven in \cite{Mark}).

\bigskip
\noindent{\sc Proposition 3.4~:}\enskip\enskip {\it Let $G$ be a homeomorphism of $\T^2$ isotopic to $L$ and $H$ the unique factor map from $G$ to $L$ that is homotopic to the identity. Then every fiber $H^{-1}(\{\overline a\})$, $\overline a\in\T^2$, is an acyclic compact set.}

\bigskip \noindent{\it Proof.}\enskip Let $g$ be a lift of $G$ to $\R^2$ and $h$ the unique lift of $H$ that defines a factor map from $g$ to $l$. The map $h$ is proper, so each set $k_a=h^{-1}(a)$ is compact. Let us prove that its complement is connected, which means that $k_a$ is acyclic. Let us begin with some remarks.

If $z$, $z'$ in $\R^2$ have the same image by $h$, then the sequence $\left( g^n(z)-g^n(z')\right)_{n\in\Z}$ is bounded by $2C_*C$, because
$$\eqalign{\Vert g^n(z)-g^n(z')\vert&\leq \Vert g^n(z)-l^n(h(z))\Vert + \Vert l^n(h(z'))-g^n(z')\Vert\cr
&\leq 2C_*C.\cr}$$

If $z$, $z'$ in $\R^2$ have not the same image by $h$, then the sequence  $\left(g^n(z)- g^n(z')\right)_{n\in\Z}$ is not bounded because
$$\eqalign{\Vert g^n(z)-g^n(z')\Vert
&\geq\Vert l^n(h(z))-l^n(h(z'))\Vert-\Vert g^n(z)-l^n(h(z))\Vert - \Vert l^n(h(z'))-g^n(z')\Vert\cr
&\geq\Vert l^n(h(z))-l^n(h(z'))\Vert-2C_*C,}$$
and because one least of the two sequences 
$$\left(\Vert l^n(h(z))-l^n(h(z')\right)\Vert)_{n\geq 0} 
\enskip \mathrm{or} \enskip\left(\Vert l^{-n}(h(z))-l^{-n}(h(z'))\Vert\right)_{n\geq 0} $$ tends to $+\infty$ when $n$ tends to $+\infty$. 

Let us fix $a\in\R^2$. To prove that $\R^2\setminus k_a$ is connected, we need to prove that every point $z'\in\R^2\setminus k_a$ belongs to the unbounded connected component of $\R^2\setminus  k_a$ or equivalently that there exists $n\in\Z$ such that $g^n(z')$ belongs to the unbounded connected component of $\R^2\setminus  g^n(k_a)$.  Choose $z\in k_a$. There exists $n\in\Z$ such that $\Vert g^n(z)-g^n(z'))\Vert>2C_*C $, which implies that $g^n(z')$ belongs to the unbounded connected component of $\R^2\setminus  g^n(k_a)$ because $g^n(k_a)=k_{l^n(a)}$.

The preimage of  $K_{\overline a}$ by the covering projection may be written $\bigsqcup_{a+\Z^2=\overline a} k_a$ and similarly the  preimage of every connected component $K'$ of $K_{\overline a}$ may be written $\bigsqcup_{\omega\in\Z^2} t_{\omega}(k')$ where $k'$ is a cellular set.  This implies that the set $K'$ itself is cellular. \hfill$\Box$

\bigskip If $G$ satifies the hypothesis of the previous proposition and if $F$ is a homeomorphism of $\T^2$ isotopic to the identity that commutes with $G$, it is easy to prove that $F$ fixes every fiber $H^{-1}(\{y\})$, $y\in\T^2$. The following result is much more surprising and is one of the fundamental results of \cite{Mark}. We will see here how to get a short proof by applying Theorem 1.1.

\bigskip
\noindent{\sc Proposition 3.5~:}\enskip\enskip {\it Let $G$ be a homeomorphism $\T^2$ isotopic to $L$, and $H$ the unique factor map from $G$ to $L$ that is homotopic to the identity. If $F$ is a homeomorphism of $\T^2$ isotopic to the identity that commutes with $G$, then $F$ fixes the connected components of the fibers $H^{-1}(\{\overline a\})$, $\overline a\in\T^2$.
}

\bigskip
\noindent{\it Proof.} \enskip  The transformation $H\circ F$ is continuous and homotopic to the identity. Moreover, one has $$H\circ F\circ G=H\circ G\circ F = L\circ H\circ F,$$ which implies that $H=H\circ F$, by uniqueness of the factor map given by Proposition 3.3. Fix a lift $g$ of $G$ to $\R^2$. Like in the proof of Proposition 3.1, the fact that $L$ has a unique fixed point implies that there exists a unique lift of $F$ to $\R^2$ that commutes with $g$. Indeed, fix a lift $f$ of $F$ to $\R^2$. Here again, there exists ${\omega}\in\Z^2$ such that $g\circ f=t_{{\omega}}\circ f\circ g$ and one must find ${\omega}'\in\Z^2$ such that $t_{{\omega}'}\circ f\circ g=g\circ t_{{\omega}'}\circ f$. As we have
$$g\circ t_{{\omega}'}\circ f=t_{l({\omega}')}\circ g\circ f=t_{l({\omega}')+{\omega}}\circ f\circ g,$$ we must solve$${\omega}'=l({\omega}')+{\omega},$$ and we have seen in the proof of Proposition 3.2 that this equation has a unique solution.

\medskip
Suppose now that $f$ is the lift of $F$ that commutes with $g$. Denote by $h$ the unique factor map from $g$ to $l$ that commutes with the integer translations. Observing that
$$h\circ f\circ g=h\circ g\circ f=l\circ h\circ f$$  and that $h\circ f$ commutes with the integer translations, one deduces that $h\circ f=h$. According to Proposition 3.4, one knows that the decomposition of $\R^2$ in fibers of $h$ is an invariant acyclic decomposition of $f$, and to Proposition 2.9, that the induced cellular decomposition ${\cal D}'=(K'_{i'})_{i'\in I'}$ is $f$-periodic. Recall that $I'$ is a topological plane by Moore's Theorem. The induced map $\varphi'$ being periodic, it has a fixed point by Brouwer's Lemma on Translation Arcs. The natural action of $\Z^2$ on $\R^2$ induces an action of $\Z^2$ on $I'$ and of course $\varphi'$ commutes with the induced  transformations. In particular, $\varphi'$ has infinitely many fixed points. We have seen in Proposition 1.6 that this implies that $\varphi'$ is the identity .\hfill$\Box$

\bigskip
\bigskip
\noindent{\bf Construction of local Anosov homeomorphisms, main result}

\bigskip

Consider the compact one punctured torus $T_0$ obtained by blowing up the origin of $\T^2$. In other words, one adds to $\T^2\setminus\{0\}$ the circle $S_0$ of directions at $0$, equipped with the natural topology. Every diffeomorphism of class $C^1$ of $\T^2$ that fixes $0$  has a natural extension to $\T^2$: the homeomorphism that acts on $S_0$ by the natural action of the derivative $DF(0)$. Denote by $L_0$ the homeomorphism obtained from our hyperbolic automorphism $L$ in that way.  The restriction $L_0\vert_{S_0}$ has four fixed points: two sinks  $\omega$, $\omega'$ separated by two sources $\sigma$, $\sigma'$. Let us denote by $I_r$, $1\leq r\leq 4$, the connected components of $S_0\setminus\{\omega,\omega',\sigma,\sigma'\}$ and choose, for every $r$, an increasing homeomorphism $h_r~:\overline I_r\to[0,1]$, where $\overline I_r$ inherits the original orientation of $S_0$. We define a homeomorphism $L'_0$ of $S_0\times [0,1]$ by setting 
$$ L'_0(z,s)= (h_r^{-1}( (1-s) h_r(L_0(z)) + s h_r(z)),s),$$for every $(z,s)\in I_r\times[0,1]$. We obtain a compact one punctured torus $T_1$ by pasting
$T_0$ and $S_0\times [0,1]$ identifying $z\in S_0$ with $(z,0)\in S_0\times [0,1]$ and a homeomorphism  $L_1$ on $T_1$ that coincides with $L_0$ on $T_0$ and with $L'_0$ on $S_0\times [0,1]$. The homeomorphism does not depend, up to conjugacy, on the above construction.

 \medskip Let $M$ be a closed surface without boundary and  $T\subset M$ a compact one punctured torus. We will call {\it local Anosov homeomorphism supported on $T$} every homeomorphism of $M$ that fixes every point outside $T$ and that is conjugate to $L_1$ when restricted to $T$. 
 
 \medskip
 The goal of this  section is to prove the following result:

\bigskip
\noindent{\sc Theorem 3.6 ~:}\enskip\enskip  {\it Let $M$ be an orientable  closed surface with no boundary, $T\subset M$ a compact one punctured torus, $\check T$ the $2$-torus obtained from $M$ by identifying $\overline {M\setminus T}$ to a point and $P~: M\to \check T$ the natural projection. Let $F$ be a homeomorphism of $M$ such that:

\smallskip 
\noindent{\bf i)}\enskip\enskip $F$ is isotopic to a homeomorphism that fixes every point of $T$;

\smallskip 
\noindent{\bf ii)}\enskip\enskip  there exists a homeomorphism that commutes with $F$ and is isotopic in $M$ to a local Anosov homeomorphism supported on $T$.

\smallskip
Then, there exists a connected component $U$ of the domain of the canonical invariant essential acyclic decomposition of $F$ such that the map $P^*: H^1(\check T,\Z)\to H^1(U,\Z)$ is injective.
}

\bigskip
\bigskip
\noindent{\bf Study of local Anosov homeomorphisms}

\bigskip
 
We will extend what we did for hyperbolic automorphisms of $\T^2$ to local Anosov homeomorphisms. We fix here a closed surface without boundary $M$ and a local Anosov homeomorphism $L^*$ supported on a compact one punctured torus $T$. To simplify the notations, we will write $T=T_1$ and $L^*\vert _T=L_1$, keeping the notations of the previous subsections. We will give the same name to the projection $P$ defined in the statement and the projection $P~: M\to\T^2$ that sends every point of  $\T^2\setminus\{0\}=T_0\setminus \partial T_0$ on itself and every point of $N=\overline {M\setminus T_0}$ on $0$. Indeed, the map $P~: M\to\T^2$ defines a natural map $\check P~: \check T\to \T^2$ such that $\check P^*: H^1(\T^2,\Z)\to H^1(\check T,\Z)$ is an isomorphism. Thus, to prove Theorem 3.6, it is sufficient to prove that there exists a connected component $U$ of the domain of the canonical invariant essential acyclic decomposition of $F$ such that the map $P^*: H^1(\T^2,\Z)\to H^1(U,\Z)$ is injective. Fix $z_0\in \T^2\setminus\{0\}$. The normal covering space $\widetilde M$ of $M$ obtained as the quotient of the universal covering space of $M$ by the kernel of  $$P_*~: \pi_1(z_0,M)\to \pi_1(z_0,\T^2)$$ has a deck transformation group $\{t^*_{\omega}\,,\, {\omega}\in\Z^2\}$ isomorphic to $\Z^2$, and $P$ can be lifted to an application $p~:\widetilde M\to \R^2$ satisfying $p\circ t^* _{\omega} =t_{\omega}\circ p$, for every ${\omega}\in\Z^2$. The inverse image $p^{-1}(\{\omega\})$ of a point ${\omega}\in\Z^2$ is a surface with boundary $N_{\omega}$ and the covering projection $\Pi~: \widetilde M\to M$ induces a homeomorphism between $N_{\omega}$ and $N$. By contrast, the inverse image of a point ${\omega}\not\in\Z^2$ is reduced to a point. The homeomorphism $L^*$ can be lifted to $l^*$, where  $l^*\circ t^*_{\omega}=t^*_{l({\omega})}\circ l^*$, and where $l^*$ coincides on a set $N_{\omega}$, ${\omega}\in\Z^2$, with the transformation $t^*_{l({\omega})-{\omega}}$. In particular, the fixed point set of $l^*$ is included in $N_0$.

\medskip
 
 The following result is the natural extension of Proposition 3.3.

\bigskip
\noindent{\sc Proposition 3.7~:}\enskip\enskip {\it Let $G$ be a homeomorphism of $M$ isotopic to $L^*$. Then, there exists a unique continuous map $H~: M\to\T^2$ homotopic to $P$, such that $H\circ G=L\circ H$, and this map is onto. }

\bigskip \noindent{\it Proof.}\enskip Let $(G_s)_{s\in[0,1]}$ be an isotopy in $M$ from $L^*$ to $G$. One can lift this isotopy on $\widetilde M$ to an isotopy $(g_s)_{s\in[0,1]}$ such that $g_0=l^*$. Of course, one has $g_s\circ t^*_{\omega}=t^*_{l({\omega})}\circ g_s$, for every ${\omega}\in \Z^2$. Observe that there exists $C>0$ such that for every $s\in[0,1] $ and every $z\in\widetilde M$, one has
$$\Vert p\circ g_s(z)-l\circ p(z)\Vert=\Vert p\circ g_s(z) - p\circ l^*(z)\Vert\leq C.$$
This implies that the sequence $(p\circ g_s^k(z))_{k\in\Z}$ is a $C$-pseudo-orbit of $l$. By Theorem 3.1 there exists a unique point $z'\in \R^2$ such that the sequence $(p\circ g_s^k(z)-l^k(z'))_{k\in\Z}$ is bounded. Moreover, this sequence is bounded by $C_*C$. We denote by $h_s(z)$  this point. Like in the proof of Proposition 3.3, the graph of the function $(s,z)\mapsto h_s(z)$ is the following closed set
$$\left\{(s,z,z')\in[0,1]\times\widetilde M\times\widetilde M\enskip\vert\enskip k\in\Z\Rightarrow\Vert p\circ g_s^k(z)-l^k(z')\Vert \leq C_*C\right\},$$
which implies that the map $(s,z)\mapsto h_t(z)$ is continuous on $[0,1]\times\widetilde M$. In this new situation, one gets

\smallskip
\noindent-\enskip\enskip  $h_0=p$;

\smallskip
\noindent-\enskip\enskip   $h_s\circ g_s=l\circ h_s$;

\smallskip
\noindent-\enskip\enskip   $h_s\circ t^*_{{\omega}}=t_{{\omega}}\circ h_s$ for every  $\omega\in\Z^2$.

\smallskip
Therefore, $(h_s)_{s\in[0,1]}$ lifts a continuous family $(H_s)_{s\in[0,1]}$ of maps from $M$ to $\T^2$, such that $H_0=P$ and such that $H_s\circ G_s=L\circ H_s$ for every $s\in[0,1]$. The map $H_0=P$ inducing a group isomorphism between $H_2(M,\Z)$ and $H_2(\T^2,\Z)$, it is the same for $H=H_1$, which implies that $H$ is onto.

 \medskip
To conclude the proof, it remains to prove the unicity of the factor map $H$, which is done like in the proof of Proposition 3.3. One must prove that every continuous map $H'~: M\to\T^2$ homotopic to $P$ such that $H'\circ G=L\circ H'$ coincides with $H$. Consider a continuous family $(H'_s)_{s\in[0,1]}$ of maps from $M$ to $\T^2$ such that $H'_0=P$ and $H'_1=H'$. Lift it to a continuous family $(h'_s)_{s\in[0,1]}$ of maps from $\widetilde M$ to $\R^2$ such that $h'_0=p$. Set $h'=h'_1$. Observe that  $h'\circ t^*_{\omega}=t_{\omega}\circ h'$ for every ${\omega}\in\Z^2$.  As a consequence of the equality $L\circ H'=H'\circ G$, one knows that there exists ${\omega}_0\in\Z^2$ such that $l\circ  h'=t_{{\omega}_0}\circ h'\circ g$. By hypothesis, $L$ has a unique fixed point and ${\rm det}(l-{\rm Id})=\pm 1$. One deduces that there exists ${\omega}_1\in\Z^2$ unique, such that $l({\omega}_1)-{\omega}_1=-{\omega}_0$. Set $h''=t_{{\omega}_1}\circ h'$ and observe that
$$l\circ h''=l\circ t_{{\omega}_1}\circ h'=t_{l({\omega}_1)}\circ l\circ h'=t_{l({\omega}_1)}\circ t_{{\omega}_0}\circ h'\circ g=t_{{\omega}_1}\circ h'\circ g=h''\circ g.$$ One deduces that $l^k\circ h''=h''\circ g^k$ for every $k\in\Z$. The map $h''-p$ being invariant by the deck transformations $t^*_{\omega}$, $\omega\in\Z^2$, it is bounded. This implies that $h''(z)$ is the unique point $z'$ such that the sequence $(p\circ g^k(z)-l^k(z'))_{k\in\Z}$ is bounded, that means $h''(z)=h(z)$. One deduces that $H'=H$. \hfill$\Box$

\medskip
Now, me must extend proposition 3.4:

\bigskip
\noindent{\sc Proposition 3.8~:}\enskip\enskip {\it We keep the same notations as in Proposition 3.7.
Then the map $h$ is proper and every set $k_a=h^{-1}(a)$, $a\in\R^2$, is a compact set whose complement is connected. Moreover, $k_a$ is acyclic if $a\not\in\Z^r$.}

\bigskip \noindent{\it Proof.}\enskip The map $p$ being proper, it is the same for $h$ because $\Vert p(z) - h(z)\Vert\leq C_*C$ for every $z\in\widetilde M$. Thus $k_a=h^{-1}(a)$ is compact and it remains to prove that its complement is connected. Like in the proof of Proposition 3.4, let us begin with some remarks.

If $z$, $z'$ in $\widetilde M$ have the same image by $h$, then the sequence $(p(g^n(z))-p( g^n(z')))_{n\in\Z}$ is bounded by $2C_*C$, because
$$\Vert p(g^n(z))-p(g^n(z'))\Vert
\leq  \Vert p(g^n(z))-l^n(h(z))\Vert + \Vert l^n(h(z'))-p(g^n(z'))\Vert
\leq 2C_*C.$$

If $z$, $z'$ in $\widetilde M$ have not the same image by $h$, then the sequence $(p(g^n(z))-p( g^n(z')))_{n\in\Z}$ is not bounded because
$$\eqalign{\Vert p(g^n(z))-p(g^n(z'))\Vert &\geq   \Vert l^n(h(z))-l^n(h(z'))\Vert -\Vert p(g^n(z))-l^n(h(z))\Vert -\Vert l^n(h(z'))-p(g^n(z'))\Vert \cr &
\geq\Vert l^n(h(z))-l^n(h(z'))\Vert -2C_*C,\cr}$$
and because one least of the two sequences 
$$\left(\Vert l^n(h(z))-l^n(h(z')\right)\Vert)_{n\geq 0} 
\enskip \mathrm{or} \enskip\left(\Vert l^{-n}(h(z))-l^{-n}(h(z'))\Vert\right)_{n\geq 0} $$ tends to $+\infty$ when $n$ tends to $+\infty$. 

For every $r>0$, there exists $R_r>r$ such that, if $z$, $z'$ in $\widetilde M$ satisfy $\Vert p(z)-p(z')\Vert>R_r$, then $z'$ is in the unbounded connected component of $\widetilde M\setminus B(z,r)$, where we set $$B(z,r)=\{z'\in\widetilde M\,\vert\, \Vert p(z)-p(z')\Vert\leq r\}.$$

Let fix $a\in\R^2$. To prove that $\widetilde M\setminus k_a$ is connected, we need to prove that every point $z'\in\widetilde M\setminus k_a$ belongs to the unbounded connected component of $\widetilde M\setminus  k_a$ or equivalently that there exists $n\in\Z$ such that $g^n(z')$ belongs to the unbounded connected component of $\widetilde M\setminus  g^n(k_a)$.  Choose $z\in k_a$. There exists $n\in\Z$ such that $\widetilde d(g^n(z), g^n(z'))>R_{2C_*C}$, which implies that $g^n(z')$ belongs to the unbounded connected component of $\widetilde M\setminus  B(g^n(z), 2C_*C)$. This implies that $g^n(z')$ belongs to the unbounded connected component of $\widetilde M\setminus  g^n(k_a)$ because $g^n(k_a)=k_{l^n(a)}\subset B(g^n(z), 2C_*C)$

\medskip
Let us prove now the second assertion and suppose that $a\not\in\Z^2$. For every $\omega\in\Z^2$, there exists $n\in\Z$ such that $\Vert l^{n}(\omega) -l^n(a)\Vert>C_*C$. Fix $z\in N_{\omega}$ and $z'\in k_a$. One has
$$\eqalign{ \Vert p(l^{*n}(z))-p(g^n(z'))\Vert
&\geq   \Vert p(l^{*n}(z))-l^n(h(z'))\Vert -\Vert l^n(h(z'))-p(g^n(z'))\Vert\cr 
&\geq  \Vert l^{n}(\omega))-l^n(a)\Vert -C_*C>0\cr }.$$ 

We have proven that $l^{*n}(N_{\omega})\cap g^n(k_a)=\emptyset$ and we deduce that $g^{-n}\circ l^{*n}(N_{\omega})\cap k_a=\emptyset$. Recall that $g^{-n}\circ l^{*n}$ is isotopic to the identity. One deduces that every surface $N_{\omega}$, $\omega\in\Z^2$, is isotopic in $\widetilde M$ to a surface that is disjoint from $k_a$. As we know that $\widetilde M\setminus k_a$ is connected, it remains to prove that $k_a$ is contained in a disk to ensure that it is acyclic. Denote by $\cal N$ the set of compact connected surfaces in $\widetilde M$ with boundary, whose interior contains $k_a$. This set is non empty because every compact connected surface with boundary, `` sufficiently large'', will belong to $\cal N$. To prove that $k_a$ is contained in a disk, it is sufficient to prove the following:

\smallskip 
\noindent -\enskip\enskip if $N\in{\cal N}$ has a positive genus $g$, then there exists $N'\in{\cal N}$ whose genus $g'$ is smaller than $g$;

\smallskip 
\noindent -\enskip\enskip if $N\in{\cal N}$ has genus $0$ and if its complement has $n\geq 2$ connected components, then there exists $N'\in{\cal N}$ with genus $0$ such that its complement has $n'<n$ connected components.

\smallskip 
\noindent -\enskip\enskip if $N\in{\cal N}$ has genus $0$ and a connected complement and if the boundary of $N$ has $m\geq 2$ components, then there exists $N'\in{\cal N}$ with genus $0$ and a connected complement such that its boundary has $m'<m$ components.

\medskip 
Let us prove the first statement. Suppose that $N\in{\cal N}$ has a positive genus $g$. One can find two simple closed curves $\Gamma$ and $\Gamma'$ in $N$ that intersect in a unique point with a tranverse intersection. The cycle $\Gamma$ is homologous in $H_1(\widetilde M,\Z)$ to a cycle $\sum_{1\leq r\leq R} \Gamma_r$, where each $\Gamma_r$ is a simple closed curve included in a surface $N_{\omega_r}$. Each curve $\Gamma_r$ is isotopic to a simple closed curve $\Gamma'_r$ disjoint from $k_a$, because $N_{\omega_r}$ is isotopic in $\widetilde M$ to a surface that is disjoint from $k_a$. There exists at least one curve $\Gamma_r$ whose algebraic intersection number with $\Gamma'$ is not zero and we have the same property for $\Gamma'_r$. Either $\Gamma'_r$ is included in the interior $N$, or there is a sub-arc $\gamma'_r$ of $\Gamma'_r$ contained in the interior of $N$ except the ends which are on the boundary of $N$. If we cut the surface $N$ along $\Gamma'_r $ in the first case, along $\gamma'_r$ in the second one, we get a surface $N'\in{\cal N}$ whose genus $g'$ is smaller than $g$.

\medskip
Let us prove now the second statement. Suppose that $N\in{\cal N}$ has genus $0$ and that its complement has $n\geq 2$ connected components. One can find a simple arc in $\widetilde M\setminus k_a$ whose ends belong to distinct connected components of $\widetilde M\setminus N$, because $\widetilde M\setminus k_a$ is connected. There exists a sub-arc $\gamma'$ of $\gamma$ that is contained in the interior or $N$ except the ends which are on the boundary of $N$ and belong  to the closure of distinct connected components of $\widetilde M\setminus N$. If we cut the surface $N$ along $\gamma'$, we get a surface $N'\in{\cal N}$ with genus $0$ such that its complement has $n'<n$ connected components.

\medskip
Let us finish with the last statement. Suppose that $N\in{\cal N}$ has genus $0$,  that its complement is connected and that the boundary of $N$ has $m\geq 2$ components. Fix such a component $\Gamma'$. One can find a simple closed curve $\Gamma$ in $\widetilde M$ that intersects $\Gamma'$ in a unique point with a tranverse intersection. Like in the first case, the cycle $\Gamma$ is homologous in $H_1(\widetilde M,\Z)$ to a cycle $\sum_{1\leq r\leq R} \Gamma'_r$, where each $\Gamma'_r$ is a simple closed curve disjoint from $k_a$. There exists at least one curve $\Gamma'_r$ whose algebraic intersection number with $\Gamma'$ is not zero. There is a sub-arc $\gamma'_r$ of $\Gamma'_r$ contained in the interior of $N$ except the ends which are on the boundary of $N$ and belong to distinct connected components of the boundary of $N$. If we cut the surface $N$ along $\gamma'_r $ we get a surface $N'\in{\cal N}$ with genus $0$ and a connected complement such that its boundary has $m'<m$ components. \hfill$\Box$

\bigskip
\noindent{\sc Corollary 3.9~:}\enskip\enskip {\it We keep the same notations as in Proposition 3.7.
Then, every set $K_{\overline a}=H^{-1}(\overline a)$, $\overline a\in\T^2\setminus\{0\}$, is acyclic.}

\bigskip \noindent{\it Proof.}\enskip Here again, the preimage of  $K_{\overline a}$ by the covering projection may be written $\bigsqcup_{a+\Z^2=\overline a} k_a$ and similarly the preimage of every connected component $K'$ of $K_{\overline a}$ may be written $\bigsqcup_{\omega\in\Z^2} t^*_{\omega}(k')$ where $k'$ is a cellular set.  This implies that the set $K'$ itself is cellular. \hfill$\Box$

\bigskip
Let us conclude this section with the proof of Theorem 3.6.

\bigskip

\noindent{\it Proof of Theorem 3.6.}\enskip \medskip
By Corollary 3.9, one knows that the domain $W$ of the canonical invariant acyclic decomposition of $F$ contains $M\setminus K_{\overline 0}$. There exists a neighborhood $V$ of  $K_{\overline 0}$ such that for every compact surface with boundary $N$ satisfying $K_{\overline 0}\subset N\subset V$, the connected components of $\pi^{-1}(N)$ are compact. We fix such a surface $N$ and want to find a connected component $V$ of $M\setminus N$ such that $P^*: H^1(\T^2,\Z)\to H^1(V,\Z)$ is injective. Write $(x_1,x_2)$ for the  coordinates on $\R^2$. The one-forms $dx_1$, $dx_2$ on $\T^2$ defines the generators $\lambda_1$, $\lambda_2$ of $H^1(\T^2,\Z)$ and the two-form $dx_1\wedge dx_2$ the generator $\lambda_1\smile\lambda_2$ of $H^2(\T^2,\Z)$. The fact that the connected components of $\pi^{-1}(N)$ are compact implies that $P^*(\lambda_1)$ and $P^*(\lambda_2)$ vanish on $H_1(N,\Z)$ and therefore define elements of $H^1(M,N,\Z)$. By duality, to each cohomology class $P^*(\lambda_i)\in H^1(M,N,\Z)$, $i\in\{1,2\}$, is associated an homology class $\alpha_i\in H_1(M\setminus N)$ (see Bredon \cite{Bre}, Corollary 8.4, for example) and to $P^*(\lambda_1)\smile P^*(\lambda_2)= P^*(\lambda_1\smile\lambda_2)\in H^2(M,N,\Z)$ is associated by duality the intersection $\alpha_1\wedge\alpha_2\in H_0(M\setminus N)$. Denote by $(V_r)_{r\in R}$ the (finite) family of connected components of $M\setminus N$ and for every $r\in R$,  write $\nu_r$ for the natural generator of $H_0(V_r,\Z)$. One may write $\alpha_1\wedge\alpha_2=\sum_{r\in R}k_r\nu_r$, where $k_r\in \Z$. One knows that $\sum_{r\in R}k_r=1$ and that $P^*: H^1(\T^2,\Z)\to H^1(V_r,\Z)$ is injective if $k_r\not=0$.\hfill$\Box$
 
 \bigskip

\noindent{\bf Remark:}\enskip The decomposition that we have constructed above depends on $G$ (and not on $F$). Thus, we have proved something stronger that what is stated in Proposition 3.6. Let $M$ be an orientable  closed surface, $T\subset M$ a compact one punctured torus and $G$ a homeomorphism isotopic to a local Anosov homeomorphism supported on $T$. Denote by $\cal F$ the group of homeomorphisms of $M$ that commute with $G$ and are isotopic to homeomorphisms that fix every point of $T$. Then, there exists a connected component $\cal U$ of the domain of the canonical invariant essential acyclic decomposition of the group $\cal F$, such that $P^*: H^1(\check T,\Z)\to H^1({\cal U},\Z)$ is injective.
 
\bigskip
\bigskip

\noindent{\large {\bf 4. Non existence of a section of the mapping class group}}

\bigskip
In this section, we fix an orientable closed surface $M$, we denote by ${\rm Homeo}_+(M)$ the group of orientation preserving homeomorphisms of $M$ and ${\rm MCG}(M)$ the Mapping Class Group of $M$ with the group structure naturally induced. We denote by $${\cal P} : {\rm Homeo}_+(M)\to {\rm MCG}(M)$$ the projection, which associates to every homeomorphism $f$ its isotopy class $[f]$. We want to prove:

\bigskip
\noindent{\sc Theorem 4.1~:}\enskip\enskip {\it If the genus $g$ of $M$ is larger than $1$, there is no morphism   $${\cal E} :  {\rm MCG}(M)\to{\rm Homeo}_+(M)$$ such that $${\cal P}\circ {\cal E}={\rm Id}_{{\rm MCG}(M)}.$$ }

\bigskip
The proof will be done by contradiction, by supposing that $\cal E$ exists.  Denote by $\Theta$ the {\it fundamental Dehn twist}, which is the homeomorphism of $\T^1\times[0,1]$ that is lifted to the universal covering space by $$\eqalign{\theta  ~: \R\times [0,1]&\to \R\times [0,1],\cr (x,y)&\mapsto (x+y,y).\cr}$$ For every simple closed curve $\beta$ that is not homotopic to the identity we choose a one to one continuous map  $H_{\beta}~: \T^1\times[0,1]\to M$, whose image will be denoted by ${\cal A}_{\beta}^*$, such that every curve $H_{\beta}(\T^1\times\{y\})$ is freely homotopic to $\beta$, then we define $F^*_{\beta}\in {\rm Homeo}_+(M)$ which is equal to the identity outside ${\cal A}_{\beta}^*$ and is conjugate to $\Theta$ by $H_{\beta}^{-1}$ on this annulus. The isotopy class $[F_{\beta}^*]$ is independent of the choice of $H_{\beta}$, of the orientation of $\beta$ and of the free homotopy class of $\beta$. We define $F_{\beta}={\cal E}([F_{\beta}^*])$.  Let us apply first what has been done in the previous section to the maps $F_{\beta}$. Let denote by $W_{\beta}$ the domain of the canonical invariant essential cellular decomposition of  $F_{\beta}$.

\bigskip
\noindent{\sc Proposition 4.2~:}\enskip\enskip {\it  There exists an annular set ${\cal A}_{\beta}$, intersection of a nested sequence of compact annuli homotopic to ${\cal A}_{\beta}^*$ such that $U_{\beta}=M\setminus {\cal A}_{\beta}$ is one of the connected component of $W_{\beta}$, the other components being essential annular open sets included in ${\cal A}_{\beta}$. }

\bigskip
 \noindent{\it Proof.}\enskip Let us begin by proving that the domain $W'_{\beta}$ of the canonical invariant essential acyclic decomposition of  $F_{\beta}$ is not empty and for that, fix a one punctured compact torus $T\subset W_{\beta}^*$. One can construct a local Anosov homeomorphism $L^*$ supported on $T$. The supports of $F_{\beta}^*$ and $L^*$ being disjoint, these two maps commute. One deduces that $F_{\beta}^*$ and $[L^*]$ commute which implies that $F_{\beta}$ and ${\cal E}(L^*)$ commute. The hypotheses of Theorem 3.6 are satisfied for $F_{\beta}$ and $T$. Thus, there exists a connected component $U$ of $W'_{\beta}$ such that the map $P^*: H^1(\check T,\Z)\to H^1(U,\Z)$ is injective, where $P~: M\to \check T$ is the natural projection on the  $2$-torus $\check T$ obtained from $M$ by identifying $\overline {M\setminus T}$ to a point. 
 
 \medskip
 
Now, let us study the topological properties of $W'_{\beta}$. According to Corollary 2.12 and to the remark that follows it, one can write $W'_{\beta}=\bigcup_{n\geq 0} Q_n$ where $(Q_n)_{n\geq 0}$ is a sequence of  compact surfaces with boundaries, such that
 $Q_{n}\subset{\rm Int }(Q_{n+1})$, and such that $Q_{n+1}\setminus{\rm Int}(Q_n)$ is a union of essential compact annuli. We will prove that every boundary circle of $Q_n$ is homotopic to $\beta$. Consider a connected component $U$ of $W'_{\beta}$. According to the remark that follows Theorem 2.10 and Cartwright-Littlewood's Fixed Point Theorem, one knows that there exists an integer $q\geq 1$ such that every end of $U$ is accumulated by fixed points of $F_{\beta}^q$. To such an end is naturally associated one boundary circle of $Q_n$ and conversely to each boundary circle is associated naturally one end of  $W'_{\beta}$. One deduces that every boundary circle of $Q_n$ is freely homotopic to its image by a power of $F_{\beta}$. Fix $n_0$ and choose a boundary circle $\Gamma$ of $Q_{n_0}$. As $F_{\beta}$ is isotopic to $F^*_{\beta}$, one deduces that there exists $q\geq 1$ such that $F^{*q}_{\beta}(\Gamma)$ is homotopic to $\Gamma$. The homeomorphism $F^*_{\beta}$ being a Dehn twist, one deduces that $\Gamma$ is homotopic to a simple closed curve $\Gamma'\in W^*_{\beta}$. It remains to prove that $\Gamma'$ is homotopic to $\beta$. If $\Gamma'$ does not separate $W^*_{\beta}$, one can find a compact one punctured torus $T'\subset W_{\beta}^*$ that contains $\Gamma'$ and do the same construction as above. As a consequence there exists a loop in $W'_{\beta}$ whose algebraic intersection number with $\Gamma'$ is non zero. If $n$ is sufficienty large, this loop will have an algebraic intersection number equal to zero with the boundary circle of $Q_n$ corresponding to the same end than $\Gamma$, which gives us a contradiction.  If $\Gamma'$ separates $W^*_{\beta}$ but is not homotopic to $\beta$, then one can find a non separating simple loop $\Gamma''\subset W^*_{\beta}$ whose geometric intersection number is non zero (that means that every loop homotopic to $\Gamma''$ meets $\Gamma'$) and then a compact one punctured torus $T'\subset W_{\beta}^*$ that contains $\Gamma''$. Observe that there is no connected $U$ component of $Q_{n_0}$ such that  $P^*: H^1(\check T',\Z)\to H^1(U,\Z)$ is injective. Here again we have a contradiction.

 \medskip
  We have proven that every boundary circle of $Q_n$ is homotopic to $\beta$. In fact we have proven more: either the domain $Q_n$ is equal to the whole surface $M$ or it contains a subsurface, union of one or two of its connected components (depending whether $\beta$ is separating or not), that is isotopic to the closure of $W^*_{\beta}$. In the first case, the domain $W'_{\beta}$ coincides with $M$. In the second case, there exists an annular set ${\cal A}_{\beta}$, intersection of a nested sequence of compact annuli homotopic to ${\cal A}_{\beta}^*$ such that $U_{\beta}=M\setminus {\cal A}_{\beta}$ is the union of one or two connected connected components of $W'_{\beta}$,  the other components being essential annular open sets included in ${\cal A}_{\beta}$. Of course, $U_{\beta}$ is uniquely defined, invariant  by $F_{\beta}$, and each of its two connected components is invariant by $F_{\beta}$ in the case where it is not connected. Define the set $U$ to be equal to $U_{\beta}$ in the case where $U_{\beta}$ is connected, and to one of its connected component in the opposite case. Denote by ${\cal D}_{\beta} = (\widetilde K_i)_{i\in I_{\beta}}$ the cellular decomposition induced on $U$ by the restriction to $U$ of the canonical essential acyclic decomposition of $F_{\beta}$. By Proposition 2.9, there exists $q\geq 1$ such that ${\cal D}_{\beta} $ is invariant by $F_{\beta}^q$. To get our result it remains to prove that $U\not=M$ and $q=1$. 

\medskip
 Write $\widetilde M$ for the universal covering space of $M$ and $\Pi~: \widetilde M \to M$ for the covering projection. Denote by $\cal T$ the deck transformation group. Define the set $W^*$ to be equal to $W^*_{\beta}$ in the case where $\beta$ does not separate $M$ and to the connected component of $W^*_{\beta}$ corresponding to $U$ in the case where $\beta$ separates $M$ and $U_{\beta}$ is not connected. Fix a connected component $\widetilde W^*$ of $\Pi^{-1}(W^*)$. It is the universal covering space of $W^*$ with a group of covering automorphisms equal to the stabilizer of  $\widetilde W^*$  in $\cal T$, that we denote by ${\cal T}'$. There exists a lift $\widetilde F_{\beta}^*$ of $F_{\beta}^*$ that fixes every point of  $\widetilde W^*$. This lift commutes with every $\tau\in{\cal T'}$. Fix an isotopy $I$ from $F^*_{\beta}$ to $F_{\beta}$, lift it to an isotopy $\widetilde I$ defined on $\widetilde M$ and starting from $\widetilde F^*_{\beta}$ and denote by $\widetilde F_{\beta}$ the other end of the isotopy. It is a lift of $F_{\beta}$ that commutes with every element of ${\cal T}'$. One can find a non trivial element $\tau\in{\cal T}'$ such that any path joining a point $\widetilde z$ to $\tau(\widetilde z)$ projects on a loop of $M$ that is non homotopic to $\beta$ whatever $\beta$ is oriented. There exists a unique connected component $\widetilde U$ of $U$ such that $\tau$ belongs to the stabilizer of $\widetilde U$. The fact that $\widetilde F_{\beta}$ commutes with $\tau$ implies that $\widetilde F_{\beta}(\widetilde U) =\widetilde U$. The decomposition ${\cal D}_{\beta} $  can be lifted  to a cellular decomposition $\widetilde{\cal D} = (\widetilde K_i)_{i\in \widetilde I}$ on $\widetilde U$. The map $\widetilde F$ acts on $\widetilde{\cal D}$ and there exists $\tau'\in{\cal T}'$ such that for every $i\in \widetilde I$, one has $\widetilde F^q_{\beta}(\widetilde K_i)=\tau'(\widetilde K_i)$ because ${\cal D}_{\beta} $  is invariant by $F_{\beta}^q$. The fact that $\widetilde F_{\beta}$ commutes with every $\tau\in{\cal T}'$ implies that $\tau'$ belongs to the center of ${\cal T}'$. But this center is trivial because $W^*$ is a two punctured surface of genus $g-1$ or one punctured surface of genus $g'\geq 1$. One deduces that $\widetilde{\cal D}$ is invariant by $\widetilde F_{\beta}^q$. We have meet such a situation before: the set $I'$ is a topological disk with a proper and free action of ${\cal T}'$, the  map $\varphi$ induced by $\widetilde F_{\beta}$ is a periodic map that commutes with the elements of  ${\cal T}'$, it must have a fixed point by Brouwer's Lemma on Translation Arcs and in fact infinitely many such fixed points. This implies that $\varphi$ is the identity. So we have $q=1$. As a consequence, one deduces that $\widetilde F_{\beta}$ commutes with the whole stabilizer of $\widetilde U$. As one knows that $\widetilde F_{\beta}$ does not commute with every element of $\cal T$, one deduces that $\widetilde U\not=\widetilde M$.\hfill$\Box$

\bigskip

 Now we will begin the proof of Theorem 4.1 and will suppose first that $g$ is even. We fix a simple closed curve $\beta_0$ that separates $M$ into two homeomorphic one punctured surfaces of genus $g/2$ and we consider an involution $I^*\in {\rm Homeo}_+(M)$ that fixes $\beta_0$ and permutes the two connected components of $M\setminus \beta_0$.  We know that the class $[I^*]\in {\rm MCG}(M)$ has order two, which implies that the homeomorphism $I={\cal E}(I^*)$ is an involution. In one of the connected components of $M\setminus \beta_0$ we construct a sequence $(\beta_j)_{1\leq j\leq g}$ of simple closed curves such that: 

\smallskip
\noindent-\enskip\enskip if $\vert j-j'\vert >1$, then $\beta_j$ are $\beta_{j'}$ disjoint;

\smallskip
\noindent-\enskip\enskip  if $\vert j-j'\vert =1$, then $\beta_j$ and $\beta_{j'}$ have a unique point of intersection and the intersection is transverse.

\medskip
On the other component of  $M\setminus \beta_0$, we construct another sequence $(\beta_j)_{-g\leq j\leq -1}$ of simple closed curves by setting $\beta_j=I^*(\beta_{-j})$.

\medskip
For every $j\in\{-g,\dots, g\}$, one can define the objects ${\cal A}_{\beta_j}^*$, $F_{\beta_j}^*$, $F_{\beta_j}$, ${\cal A}_{\beta_j}$, $W_{\beta_j}$, $U_{\beta_j}$, that we will write ${\cal A}_{j}^*$, $F_{j}^*$, $F_{j}$, ${\cal A}_{j}$, $W_{_j}$, $U_{_j}$ respectively to simplify the notations. 
In the case where $\vert j-j'\vert>1$, one can choose the annuli ${\cal A}_j^*$ and ${\cal A}_{j'}^*$ to be disjoint. This implies that $F_j^*$ and $F_{j'}^*$ commute. Therefore, one concludes that $[F_j^*]$ and $[F_{j'}^*]$ commute in the Mapping Class Group and that the homeomorphisms $F_j$ and $F_{j'}$ commute. For similar reasons, one knows that $F_0$ commutes with every $F_j$.  Because of the equality $\beta_j=I^*(\beta_{-j})$, one can suppose that $F_j^*$ is conjugate to $F_{-j}^*$ by $I^*$, which implies that $F_j$ is conjugate to $F_{-j}$ by $I$ and consequently that $F_0$ commutes with $I$. There is a more subtle relation, that can be read only on the isotopy classes (seeÊ \cite{FaMarg}). One has:
$$([F_1^*]\dots  [F_g^*])^{2g+2}=[F_0^*]=([F_{-1}^*]\dots [F_{-g}^*])^{2g+2},$$
which implies that
$$(F_1\circ \dots \circ F_g)^{2g+2}=F_0=(F_{-1}\circ \dots \circ F_{-g})^{2g+2}.$$

\medskip
According to Proposition 4.2, the complement $U_j$ of ${\cal A}_j$ is a two punctured surface of genus $g-1$ if $j\not=0$ and the union of two one punctured surfaces of genus $g/2$ if $j=0$. We will denote by ${\cal D}^j=\left(U_j, (K^j_i)_{i\in I_j}\right )$ the canonical invariant essential cellular decomposition of  $F_j$, restricted to $U_j$. The fact that $F_0$ commutes with $F_j$ implies that ${\cal A}_j$ is invariant by $F_0$ and that $F_0$ acts naturally on ${\cal D}^j$. For similar reasons, $I$ fixes ${\cal A}_0$, permutes the two connected components of $U_0$ and the decompositions induced by ${\cal D}^0$ on each component.

\bigskip

The following result will be very useful:

\bigskip
\noindent{\sc Lemma 4.3~:}\enskip\enskip {\it There exists a fundamental system of neighborhoods of ${\cal A}_0$ made of open annuli $A$, satisfying:

\smallskip
\noindent-\enskip\enskip   $A$ is invariant by $F_0$~;

\smallskip
\noindent-\enskip\enskip  the complement of $A$ in $M$ is  ${\cal D}^0$-saturated;

\smallskip
\noindent-\enskip\enskip  the complement of ${\cal A}_0$ in $A$ is the union of two open annuli that are  ${\cal D}^0$-saturated.}

\bigskip
 \noindent{\it Proof.}\enskip  Let $A'$ be an open annulus containing $ {\cal A}_0$. The saturation ${\cal D}^0(M\setminus A')$ is the union of two disjoint connected closed sets  invariant by $F_0$.  One deduces that the connected component $A$ of the complement of  ${\cal D}^0(M\setminus A')$ that contains $ {\cal A}_0$, is an open annulus that is invariant par $F_0$ and that $A\setminus  {\cal A}_0$ is the union of two  ${\cal D}^0$-saturated open annuli. \hfill$\Box$

\bigskip
\bigskip
\noindent{\bf Rotation number, definition of $X_{0,0}$}

\bigskip
We will use the notion of {\it rotation number} in the proof of Theorem 4.1, in fact a stronger notion that the classical one. Consider the annulus $\A=\T^1\times\R$, and denote by $$\eqalign{\Pi~:\R^2&\to \A\cr (x,y)&\mapsto(x+\Z,y)\cr}$$ the universal covering projection, $$\eqalign{\tau~:\R^2&\to \R^2\cr (x,y)&\mapsto(x+1,y)\cr}$$ the generating deck transformation, and $$\eqalign{p_1~:\R^2&\to \R\cr (x,y)&\mapsto x\cr}$$ the first projection. If $F$ is a homeomorphism of $\A$  isotopic to the identity and $f$ a lift of $F$ to $\R^2$, we will say that $z\in\A$ has a rotation number $\rho_f(z)\in\R$ if

\smallskip
\noindent{\bf i)}\enskip\enskip   the positive orbit of $z$ is relatively compact;

\smallskip
\noindent{\bf ii)}\enskip\enskip  the sequence of general term
$$ p_1(f^n(\widetilde z))-p_1(\widetilde z)-n\rho_f(z)$$ is bounded if $\widetilde z\in\Pi^{-1}(\{z\}$.

\medskip
Observe that the previous sequence does not depend on the choice of the point $\widetilde z\in\Pi^{-1}(\{z\}$, that the existence of $\rho_f(z)$ does not depend on the choice of the lift $f$ and that another choice $f'$ will give a rotation number $ \rho_{f'}(z)$ such that $ \rho_{f'}(z) -\rho_f(z)\in\Z$. Let $H$ be a homeomorphism of $\A$ that induces the identity on the first homology group $H_1(\A,\Z)$ and $h$ a lift of $H$ to $\R^2$. A consequence of {\bf i)} is the fact that the rotation number $\rho_{h\circ f\circ h^{-1}}(h(z))$ exists if it is the case for $\rho_f(z)$ and that $\rho_{h\circ f\circ h^{-1}}(h(z))=\rho_f(z)$. Consequently, one can define rotation numbers in the case of  an abstract open annulus $A$ (which means a surface homeomorphic to $\A$) as soon as we fix a generator $[\gamma]$  of $H_1(A,\Z)$: if $F$ is a homeomorphism of $A$ isotopic to the identity and $f$ a lift of $F$ to the universal covering space $\widetilde {A}$ of $A$, we will say that $z\in A$ has a rotation number $\rho$ if $H(z)$ has a rotation number $\rho$, for the lift $h\circ f\circ h^{-1} $ of $H\circ F\circ H^{-1}$, where $H~:A\to\A$ is any homeomorphism whose action in homology sends $\gamma$ on the generator of  $H_1(\A,\Z)$ naturally defined by $\tau$ and  $h~:\widetilde A\to \R^2$ any lift of  $H$.

\bigskip

The technical result below will be useful:

\bigskip
\noindent{\sc Lemma 4.4~:}\enskip\enskip {\it Let $A$ be an abstract open annulus, $F$ a homeomorphism of $A$ isotopic to the identity and $f$ a lift of $F$ to the universal covering space $\widetilde {A}$ of $A$. Let $(\gamma_r)_{r\in R}$ be  a finite family of simple arcs joining the two ends of $A$ and $K\subset A$ a compact set. Let  $X\subset A$ be a connected subset such that:

\smallskip
\noindent{\bf i)}\enskip\enskip  for every $n\geq 0$, one has $F^n(X)\subset K$;

\smallskip
\noindent{\bf ii)}\enskip\enskip for every $n\geq 0$, there exists $r_n\in R$ such that $F^n(X)\cap \gamma_{r_n}=\emptyset$;

\smallskip
\noindent{\bf iii)}\enskip\enskip  there exists a point $z_0\in X$  that has a rotation number.
\medskip

\noindent Then, every point $z\in X$ has a rotation number and $\rho_f(z)=\rho_f(z_0)$. Moreover, if $X$ is invariant by $F$, then $\rho_f(z_0)$ is an integer.}

\bigskip
 \noindent{\it Proof.}\enskip Of course, one can suppose than $A=\A$. According to {\bf i)}, the positive orbit of every point $z\in X$ is relatively compact: the first condition about existence of rotation number is satisfied. Let us prove the second one. Observe that for every $r\in R$, the decomposition of $\Pi^{-1}(\A\setminus\gamma_r)$ in connected components may be written
 $$  \Pi^{-1}(\A\setminus\gamma_r)=\bigsqcup_{k\in\Z} W^k_r,$$ where $W^k_r=\tau^k(W^0_r)$, and that $\Pi$ induces a homeomorphism between $W^k_r$ and $\A\setminus\gamma_r$.

The set $K$ being compact, every set $ W^k_r\cap \Pi^{-1}(K)$ is relatively compact in $\R^2$. This implies that there exists $C>0$ such that for every $r\in R$, for every $k\in \Z$, for every $\widetilde z, \widetilde z'\in  
W^k_r\cap \Pi^{-1}(K)$, one has
$$-C\leq p_1(\widetilde z)-p_1(\widetilde z')\leq C.$$
By hypothesis, $X$ is connected and does not meet $\gamma_{r_0}$, so one has
$$  \Pi^{-1}(X)=\bigsqcup_{k\in\Z} \widetilde X ^k,$$ where $\widetilde X^k\subset W^k_{r_0}$. One knows that $\widetilde X^k=\tau^k(\widetilde X^0)$, for every $k\in\Z$, and that $\Pi$ induces a homeomorphism between $\widetilde X^k$ and $X$.  
By hypothesis, $F^n(X)$ does not meet $\gamma_{r_n}$, so there exists $k_n\in \Z$ such that $f^n(\widetilde X^0)\subset W^{k_n}_{r_n}$.

Fix $z\in X$ and denote by $\widetilde z$ and  $\widetilde z_0$ the lifts of $z$ and $z_0$ that belong to $\widetilde X^0$ respectively. By hypothesis, we know that there exists $C'>0$ such that for every $n\geq 0$, one has
$$ -C'\leq  p_1(f^n(\widetilde z_0))-p_1(\widetilde z_0)-n\rho_f(z_0)\leq C',$$ and that
$$-C\leq p_1(f^n(\widetilde z))-p_1(f^n(\widetilde z'))\leq C.$$We deduce that for every $n\geq 0$, one has
$$ -C'-2C\leq  p_1(f^n(\widetilde z))-p_1(\widetilde z)-n\rho_f(z_0)\leq C'+2C.$$
Thus, $z$ has a rotation number and $\rho_f(z)=\rho_f(z_0)$.

\medskip
Suppose moreover than $X$ is invariant by $F$. The sequence $(r_n)_{n\geq 0} $ may be chosen to be constant equal to $r_0$. In this case, the image of $\widetilde X_0$ by $f$ is a translated $\widetilde X_k$. One deduces that $f^n(\widetilde X_0) =\widetilde X_{nk}=\tau^{nk}(\widetilde X_0)$,  for every $n\geq 0$. The estimations above imply that every point $z\in X$ has a rotation number equal to $k$.  \hfill$\Box$

\bigskip
Let $A$ be an annulus given by Lemma 4.3. Observe than every orbit of  $F_0\vert_A$ is relatively compact in $A$. By fixing an orientation of $\beta_0$, one gets a natural orientation of every essential simple closed curve of $A$, that means a generator of $H_1(A,\Z)$.  One can say that one of the connected component of $A\setminus  {\cal A}_0$  is   ``on the right''  of ${\cal A}_0$ and the other one ``on the left''. Let $f$  be a lift of $F_0\vert_{A}$ to the universal covering space $\widetilde A$. The homeomorphism $F_0$ being isotopic to the Dehn twist $F_0^*$, there exists an integer $k$ such that every point in the right connected component of $A\setminus  {\cal A}_0$ has a rotation number equal to $k$ while every point in the left connected component has a rotation number equal to $k+1$. We will say that $z\in  {\cal A}_0$ has a rotation number $\rho(z)\in\R$, if the rotation number of $z$ for $f$ exists and is equal to $k+\rho(z)$. Of course, the existence and the value of  $\rho(z)$ depends neither on the choice of $f$, nor on the choice of $A$. The existence of $\rho(z)$ does not depend on the choice of the orientation of $\beta_0$ but the other choice will permute $\rho(z)$ in $1-\rho(z)$. In particular, the set $ {\cal X}_{0,0}$ of points $z\in {\cal A}_0$ whose rotation number is $1/2$, is independent of the choice of the orientation.

\bigskip
\bigskip
\noindent{\bf Properties of $ {\cal X}_{0,0}$}

\bigskip
\noindent{\sc Lemma 4.5~:}\enskip\enskip {\it The set $ {\cal X}_{0,0}$ is non empty.}

\bigskip \noindent{\it Proof.}\enskip The involution $I^*$ preserves the orientation and lets invariant $\beta_0$. This implies that its restriction to $\beta_0$ reverses the orientation, that $I^*$ has two fixed points and that these fixed points belong to $\beta_0$. Consequently, by Proposition 1.6, the involution  $I$ being isotopic to $I^*$ has exactly two fixed points. These points belong to ${\cal A}_0$ because $I$ permutes the two components of $U_0$. The homeomorphism $F_0$, commuting  with $I$, lets invariant the fixed point set of  $I$. So, either $F_0$ fixes the two fixed points, or it permutes them. In each case a fixed point of $I$ has a rotation number. In the first case, it will be an integer;  in the second case, the two fixed points of $I$ have the same rotation number (for a given orientation of  $\beta_0$). The fact that $I^*$ reverses the orientation of $\beta_0$ implies that $I$ also ``reverses'' this orientation. As $I$ commutes with $F_0$, one deduces that if $z\in {\cal A}_0$ has a rotation number $\rho(z)$, then $I(z)$ also has a rotation number and one has $\rho(I(z))=1-\rho(z)$.  
In the case where $z$ is a fixed point of  $I$, one gets $\rho(z)=1/2$. Consequently, the fixed points of $I$ cannnot be fixed by $F_0$, they are permuted by $F_0$ and belong to  ${\cal X}_{0,0}$.  \hfill$\Box$

\bigskip
\noindent{\sc Lemma 4.6~:}\enskip\enskip {\it  For every $j\in \{-g,\dots,g\}\setminus\{0\}$, one has ${\cal X}_{0,0}\subset U_j$. }

\bigskip \noindent{\it Proof.}\enskip Fix $j\in \{-g,\dots,g\}\setminus\{0\}$. The set $ {\cal A}_j$ being annular does not contain connected components of $U_0$, which implies that there exists a simple arc $\gamma$ in $U_j$ whose ends are in different connected components of  $U_0$. Let us begin by choosing an annulus $A$, according to Lemma 4.3, that does not contain the ends of $\gamma$. This implies that there is a sub-arc $\gamma'$ of  $\gamma$ that joins the two ends of $A$.  Now let us choose another annulus $A'$ given by Lemma 4.3, that is included in $A$ and relatively compact in $A$. Suppose that ${\cal A}_j\cap {\cal X}_{0,0}\not=\emptyset$ and fix $z\in {\cal A}_j\cap {\cal X}_{0,0}$. The set ${\cal A}_j$ being connected and not contained in $ {\cal A}_0$, the connected component $X$ of ${\cal A}_j\cap A'$ that contains $z$ is not included in ${\cal A}_0$ (see for example Hocking-Young \cite{HoYou}, Theorem 2-16). Each set $F_0^n(X)$, $n\geq 0$, is included in the compact set $A'$ and is disjoint from $\gamma'$ because included in $ {\cal A}_0$. One can apply Lemma 4.4. If $f$ is a lift of  $F_0\vert_A$ to the universal covering space of $\widetilde A$, there exists $k$ such that every point in $X$ has a rotation number equal to $k+1/2$. But this set contains points of $A\setminus{\cal A}_0$ and those points have  a rotation number equal to $k$ or  $k+1$. We have found a contradiction.\hfill$\Box$

\bigskip
\noindent{\sc Lemme 4.7~:}\enskip\enskip {\it If $X\subset {\cal X}_{0,0}$ is closed and invariant by $F_0$, then for every $j\in \{-g,\dots,g\}\setminus\{0\}$, the set  ${\cal D}^j(X)$ is closed, invariant by $F_0$ and included in ${\cal X}_{0,0}$ }

\bigskip \noindent{\it Proof.}\enskip Suppose that $X\subset {\cal X}_{0,0}$ is closed and invariant by $F_0$. Fix $j\in \{-g,\dots,g\}\setminus\{0\}$. For every $z\in X$, the set  ${\cal D}^j(z)$ is cellular, so there exists a simple arc $\gamma_z$ in $M$ whose ends are in distinct connected components of $U_0$, such that $\gamma_z\cap {\cal D}^j(z)=\emptyset$. The decomposition being upper semi-continuous, there exists a neighborhood $O_z\subset U_j$ of $z$ such that for every $z'\in O_z$, one has $\gamma_z\cap {\cal D}^j(z')=\emptyset$. The set $X$ being compact, it can be covered by finitely many $O_z$. Therefore, there exists a finite family $(\gamma_r)_{r\in R}$ of arcs of $M$ whose ends are in different components of $U_0$, and such that for every $z\in X$, there exists $r\in R$ with  $\gamma_r\cap {\cal D}^j(z)=\emptyset$.  

Like in Lemma 4.5, we begin by choosing an open annulus $A$, given by Lemma 4.3, such that for every $r\in R$, there exists a sub-arc $\gamma'_r$ of $\gamma_r$ that joins the two ends of  $A$. We choose a second open annulus $A'$ given by Lemma 4.3, contained in $A$ and relatively compact in $A$. Let us begin by proving that for every $z\in X$, one has ${\cal D}^j(z)\subset{\cal A}_0$. If this is not the case, there exists a connected subset $Y\subset  {\cal D}^j(z)\cap A'$ that contains $z$ and that is not included in ${\cal A}_0$. For every  $n\geq 0$, one has $F_0^n(Y)\subset {\cal D}^j(F_0^n(z))$. So, there exists $r\in R$ such that $F_0^n(Y)\cap \gamma_r=\emptyset$. Applying Lemma 4.4, we will find a contradiction like in Lemma 4.4: if $f$ is a lift of  $F_0\vert_A$ to the universal covering space $\widetilde A$, there exists $k$ such that every point in $X$ has a rotation number equal to $k+1/2$, which is untrue because $X$ contains points of $A\setminus{\cal A}_0$ and such points have  a rotation number equal to $k$ or  $k+1$. After having proven that ${\cal D}^j(z)\subset{\cal A}_0$, the same arguments permit us, applying again Lemma 4.4, to prove finally that  ${\cal D}^j(z)\subset{\cal X}_{0,0}$. To conclude it remains to say that ${\cal D}^j(X)$ is closed because the decomposition is upper semi-continuous, and invariant by $F_0$ because $X$ is invariant by $F_0$ and $F_0$ acts on ${\cal D}^j$. \hfill$\Box$

\bigskip
Let us end with this last result.

\bigskip
\noindent{\sc Lemma 4.8~:}\enskip\enskip {\it Every closed and connected set  $X\subset {\cal X}_{0,0}$ invariant by $F_0$ separates the two connected components of $U_0$. }

\bigskip \noindent{\it Proof.}\enskip If $X\subset {\cal X}_{0,0}$ is a closed and connected set invariant by  $F_0$ that does not separate the two connected components of $U_0$, one can find 
a simple arc $\gamma$ in $M$ disjoint from $X$ whose ends are in distinct connected components of $U_0$. Using Lemma 4.3, one chooses an annulus $A$ that does not contain the ends of $\gamma$ and then a sub-arc  $\gamma'$ of $\gamma$ that joins the two ends of $A$.  If $f$ is a lift of  $F_0\vert_A$ to the universal covering space $\widetilde A$, one knows by Lemma 4.4 that every point $z\in X$ has a rotation number in $\Z$, which contradicts the hypothesis $X\subset {\cal X}_{0,0}$.\hfill$\Box$

\bigskip
\bigskip
\noindent{\bf Definition of $X_{\infty}$}

\bigskip
Let us consider the alphabet $\Xi=\{-g,\dots,g\}\setminus\{0\}$ and a word $\overline j =(j_k)_{k\geq 0}\in \Xi^{\N}$ in  which every letter appears infinitely many often. By Lemma 4.7, one can define a non decreasing sequence $(X_k)_{k\geq 0}$ of $F_0$-invariant closed sets included in ${\cal X}_{0,0}$ by the following induction relation$$X_0={\rm Fix}(I), \enskip X_{k+1}={\cal D}^{j_k}(X_k)\enskip \mathrm{if} \enskip k\geq 0.$$
The sequence $(X_k)_{k\geq 0}$ converges for the Hausdorff topology to $X_{\infty}=\overline{\bigcup_{k\geq 0} X_k}$ (note that $\bigcup_{k\geq 0} X_k$ is independant of the word $ \overline j$ because every letter appears infinitely many often). The fact that every letter $j\in\Xi$ appears infinitely many often tells us that there exists a sub-sequence $(X_{k_l})_{l\geq 0}$ of ${\cal D}^j$-saturated sets, which of course are invariant by $F_j$. One deduces that $X_{\infty}$ is invariant by $F_j$ for every $j\in\{-g,\dots,g\}$ (but not necessarily ${\cal D}^j$-saturated). Using the relation $$(F_1\circ \dots \circ F_g)^{2g+2}=F_0=(F_{-1}\circ \dots \circ F_{-g})^{2g+2},$$ and the fact that $F_0$ permutes
the two points of $X_0$, one deduces that there exists $k_*\geq 0$ such that $X_k$ has two connected components if  $k<k_*$ and is connected if $k\geq k_*$. To ensure the connectedness of $X_k$, one must wait until every letter $j>0$ has appeared at least $2g+2$ times at time $k$. By Lemma 4.8, one knows that $X_k$ separates the two connected components of  $U_0$, if $k\geq k_*$.

\bigskip
\bigskip
\noindent{\bf Prime end theory}

\bigskip
Now we fix a connected component $V$ of $U_0$. Write $V_k$, $k_*\leq k\leq \infty$, for the connected component of $M\setminus X_k$ that contains $V$. Like $V$, all these sets  are one punctured surfaces of genus $g/2$. The complement of each set $V_k$, $k_*\leq k\leq \infty$, is not reduced to a point. This implies that there is a natural compactification of  $V_k$ obtained by adding {\it the circle of prime ends} (for details, see Mather \cite{Math} ). We will recall here some of the properties of the prime end theory that we will use in the proof. Let $W\subset M$ be a one punctured surface of genus $g/2$ containing $V$ and $G$ an orientation preserving homeomorphism of $M$ which lets $W$ invariant. Let us call  {\it access arc} every arc $\gamma~: [0,1[\to W$ that has a limit $z\in{\rm Fr}(W)$ at $t=1$ and say in that case that $z$ is  {\it accessible}. There exists a natural compactification of $W$ obtained by adding a circle $S^1$, that satisfies the following:

\smallskip
\noindent{\bf i)}\enskip\enskip  every access arc  $\gamma~: [0,1[\to W$ has a limit $\zeta\in S^1$ at $t=1$ if considered in the compactification $W\sqcup S^1$;

\smallskip
\noindent{\bf ii)}\enskip\enskip  two access arcs endind up at two different points of ${\rm Fr}(W)$ end up at two differents points of $S^1$;

\smallskip
\noindent{\bf iii)}\enskip\enskip  the set of accessible points is dense in ${\rm Fr}(W)$;

\smallskip
\noindent{\bf iv)}\enskip\enskip  the sets of points of $S^1$ that are limit of access arcs is dense in $S^1$;

\smallskip
\noindent{\bf v)}\enskip\enskip  the homeomorphism $G\vert_W$ can be extended to a homeomorphism  $\overline {G\vert_W}$ of $W\sqcup S^1$.

\medskip 
The condition {\bf v)} permits us to define the {\it prime end rotation number}  $\rho_{\rm pe}(G\vert_W)\in\T^1$ which is the rotation number of the homeomorphism induced on the circle of prime ends.

\medskip
We will prove the following:

\bigskip
\noindent{\sc Lemma 4.9~:}\enskip\enskip {\it For every $k\geq k_*$, one has $\rho_{\rm pe}(F_0\vert_{V_k})=1/2+\Z$ and  $\rho_{\rm pe}(F_{j_k}\vert_{V_k})=0+\Z$.  }

\bigskip \noindent{\it Proof.}\enskip The proof of the equality $\rho_{\rm pe}(F_{j_k}\vert_{V_k})=0+\Z$ is very simple. Indeed,  $X_k$ is a ${\cal D}^{j_k}$-saturated compact subset of $U_{j_k}$. So, its has a fundamenal system of ${\cal D}^{j_k}$-saturated neighborhoods. In particular, every neighborhood of $S^1$ in $V_k\sqcup S^1$ contains an element of ${\cal D}^{j_k}$ and so contains a fixed point of $F_{j_k}$ by Cartwright-Littlewood Fixed Point Theorem. This proves that the extension of $F_{j_k}$ to $V_k\sqcup S^1$ has a fixed point on $S^1$ and consequently that $\rho_{\rm pe}(F_{j_k}\vert_{V_k})=0+\Z$. Of course, the previous argument is not valid to prove that  $\rho_{\rm pe}(F_0\vert_{V_k})=1/2+\Z$. However, a recent result of  S. Matsumoto permits to conclude (see \cite{Mats}).  This result, formulated in a more general framework, asserts that $\rho_{\rm pe}(F_0\vert_{V_k})$ is the rotation number of a Borel invariant probability measure supported on $X_k$, which in our case will necessarily be $1/2+\Z$. More precisely, let us choose an annulus $A$, given by Lemma  4.3, orient the curve $\beta_0$ in such a way that $V$ is on the right of $\gamma_0$ then consider the lift $f_0$ of $F_0\vert_A$ to the universal covering space $\widetilde A$, such that all the points in the right complement $V\cap A$ of $A\setminus {\cal A}_0$ have a rotation number equal to $0$. In this case, the rotation number of every point of $X_k$ is equal to $1/2$.  The restriction $f_0\vert_{A\cap V_k}$ has a natural extension to the universal lift of $(A\cap V_k)\sqcup S^1$ and determines a lift of $\overline{F_0\vert _{A\cap V_k}}$. It induces on the boundary line a real rotation number  $\rho_{\rm pe}(f_0\vert_{A\cap V_k})$  such that $\rho_{\rm pe}(F_{0}\vert_{V_k})=\rho_{\rm pe}(f_0\vert_{A\cap V_k})+\Z$. Choose a homeomorphism $H~:A\to\A$ whose action in homology sends $\beta_0$ on the generator of  $H_1(\A,\Z)$ naturally defined by $T$ and $h~:\widetilde A\to \R^2$ a lift of  $H$.  Matsumoto's Theorem asserts that there is a Borel probability measure $\mu$ supported on $H(X_k)$ such that
$$\rho_{\rm pe}(f_0\vert_{A\cap V_k})=\int_{\A} p_1\circ h\circ f_0\circ h^{-1} (z)-p_1( z) \,d\mu.$$
However,  the Birkhoff means of the function $p_1\circ h\circ f_0\circ h^{-1} -p_1$ converge to $1/2$ on $H(X_k)$ because every point of $X_k$ has a rotation number equal to $1/2$. This implies that $\rho_{\rm pe}(f_0\vert_{A\cap V_k})=1/2$.\hfill$\Box$

\bigskip
\noindent{\sc Lemma 4.10~:}\enskip\enskip {\it One has $\rho_{\rm pe}(F_0\vert_{V_{\infty}})=1/2+\Z$ and for every $j\in\{-g\dots,g\}\setminus\{0\}$, one has $\rho_{\rm pe}(F_{j}\vert_{V_{\infty}})=0+\Z$. }

\bigskip \noindent{\it Proof.}\enskip We will prove the first equality. The second one can be proven similarly with a slight difference due to the fact that one cannot suppose the existence of a fundamental system of $F_j$-invariant neighborhoods of ${\cal A}_0$. This lemma follows from classical results about the continuity of prime end rotation numbers (for details see \cite{L}. We keep the notations of the previous lemma and we denote by $\varphi_{[k]}$, $K\leq k\leq \infty$, the homeomorphism of the boundary line defined by the natural extension  of  $f_0\vert_{A\cap V_k}$. Matsumoto's Theorem asserts that the real rotation number of $\varphi_{[k]}$ is equal to $1/2$ if $k<\infty$ and we want to prove the same equality for $\varphi _{[\infty]}$. Suppose that this is not the case, for example suppose that it is larger. In this case, there exists $N\geq 1$ such that for every $\widetilde \zeta\in\R$, one has 
$$(\varphi_{[\infty]})^{2N}(\widetilde \zeta)>\widetilde \zeta+N+2.$$  Fix an access arc $\gamma_0~: [0,1[\to 
A\cap V_{\infty}$ ending up at $z_0\in X_{\infty}$. One constructs easily another access arc $\gamma_1~: [0,1[\to 
A\cap V_{\infty}$ 
ending up at $z_1\in X_{\infty}$ such that $F_0^n(z_1)\not=z_0$, for every $n\in\{0,\dots,2N\}$, because one can avoid a finite set. Taking a sub-arc of $\gamma_1$ if necessary, one can suppose that $F_0^n(\gamma_1)\cap \gamma_0=\emptyset$, for every $n\in\{0,\dots,2N\}$. We can say more: there exists an open disk $D_0$ containing $z_0$ and an open disk  $D_1$ containing $z_1$ such that $F_0^n(\gamma_1\cup D_1)\cap( \gamma_0\cup D_0)=\emptyset$, for every $n\in\{0,\dots,2N\}$. Choose a lift $\widetilde \gamma_0$ of $\gamma_0$. There exists a unique lift $\widetilde \gamma_1$ of $\gamma_1$ to the universal covering space of $A\cap V_{\infty}$ whose limit $\widetilde \zeta_1$ 
in $\R$ satisfies $$\widetilde \zeta_0-1<\widetilde \zeta_1<\widetilde \zeta_0.$$ 
One has $$\widetilde \zeta_0+N+1<(\varphi _{[\infty]})^{2N}(\widetilde \zeta_1).$$

If $k$ is large enough, the set $X_k$ meets the two disks $D_0$ et $D_1$. One can extend the arcs $\gamma_0$ and $\gamma_1$ in $D_0$ and $D_1$ respectively, to construct access arcs $\gamma'_0$ and $\gamma'_1$ of $X_k$. If we denote by $\widetilde \gamma'_0$ and $\widetilde \gamma'_1$ the respective lifts of $\gamma'_0$ and $\gamma'_1$ that contain $\widetilde \gamma_0$ and $\widetilde \gamma_1$ and if  we denote by $\widetilde \zeta'_0\in\R$ and $\widetilde \zeta'_1\in\R$ the limits in $\R$, one must have
$$\widetilde \zeta'_1+N+1<\widetilde \zeta'_0+N+1<(\varphi _{[k]})^{2N}(\widetilde \zeta'_1),$$
which contradicts the  fact that the rotation number of  $\varphi _{[k]}$ is  $1/2$.\hfill$\Box$

\bigskip
\noindent{\bf The contradiction}

\bigskip

 \noindent{\it Proof of Theorem 4.1}\enskip It remains to prove that Lemma 4.10 contradicts our algebraic hypothesis. For every $j\in\{-g\dots,g\}\setminus\{0\}$, denote by $\Phi_j$ the homeomorphism of $S^1$ naturally defined by the restriction of $\overline{F_j\vert_{V_{\infty}}}$ to $S^1$. The family $(\Phi_j)_{-g\leq j\leq g}$ inherits all the algebraic relations of the family $(F_j)_{-g\leq j\leq g}$. In particular:

 \smallskip
 \noindent-\enskip $\Phi_0$ commutes with $\Phi_j$;

 \smallskip
 \noindent-\enskip $\Phi_j$ and $\Phi_{j'}$ commute if $\vert j-j'\vert >1$; 
 
 \smallskip
 \noindent-\enskip  we have

$$(\Phi_1\circ \dots \circ \Phi_g)^{2g+2}=\Phi_0=(\Phi_{-1}\circ \dots \circ\Phi_{-g})^{2g+2}.$$
By Lemma 4.10, we know that each $\Phi_j$, $j\not=0$, has a fixed point but that $\Phi_0$ has none. The previous equality tells us that the $\Phi_j$, $j>0$, have no common fixed point and that this is the same for the  $\Phi_{j'}$, $j'<0$. Therefore there exists $j_0\in\{1,\dots,g-1\}$ such that~:

\smallskip
\noindent-\enskip\enskip  the $\Phi_j$, $1\leq j\leq j_0$, have a common fixed point,

\smallskip
\noindent-\enskip\enskip the $\Phi_j$, $1\leq j\leq j_0+1$, have no common fixed point

\medskip
The sets $\cap_{1\leq j\leq j_0}{\rm Fix}(\Phi_j)$ and  ${\rm Fix}(\Phi_{j_0+1})$ are disjoint non empty closed subsets of $S^1$. This implies that they are finitely many connected components of $S^1\setminus\cap_{1\leq j\leq j_0}{\rm Fix}(\Phi_j)$ that meet ${\rm Fix}(\Phi_{j_0+1})$. Denote by $I_s$, $s\in S$, the open intervals of $S^1$ that we have obtained. We know that the $\Phi_{j'}$, $j'<0$, commute with the $\Phi_j$, $j>0$. So, the sets $\cap_{1\leq j\leq j_0}{\rm Fix}(\Phi_j)$ et  ${\rm Fix}(\Phi_{j_0+1})$  are invariant by each $\Phi_{j'}$, $j'<0$. So, each $\Phi_{j'}$, $j'<0$, permutes the connected components of
 $I_s$, $s\in S$. But the rotation number of $ \Phi_{j'}$ being equal to $0$, every $\Phi_{j'}$, $j'<0$, fixes all the $I_s$, $s\in S$, and so fixes their ends. Consequently, the  $\Phi_{j'}$, $j'<0$, have a common fixed point. We have found our contradiction. \hfill$\Box$
 
 \bigskip
 \noindent{\bf Remark:} \enskip Observe that the proof above can be extended obviously to the more general situation examined in Proposition E.

\bigskip
\bigskip
\noindent{\bf The odd case}

\bigskip

The proof is very similar in the case where the genus is odd. We consider two disjoint simple closed curves $\beta_{0-}$ and $\beta_{0+}$ whose union separates $M$ in a pair of two punctured surfaces of genus $(g-1)/2$, an involution $I^*\in {\rm Homeo}_+(M)$ that fixes $\beta_{0-}$ and $\beta_{0+}$ and permutes the two connected components of $M\setminus (\beta_{0-}\cup\beta_{0+})$ and an involution $J^*\in {\rm Homeo}_+(M)$ that commutes with $I^*$ and permutes $\beta_{0-}$ and $\beta_{0+}$. Here again, in one of the connected component of $M\setminus(\beta_{0-}\cup\beta_{0+})$, we construct a sequence of simple closed curves  $(\beta_j)_{1\leq j\leq g}$ satisfying: 

\smallskip
\noindent-\enskip\enskip $\beta_j$ and $\beta_{j'}$ are invariant by $J^*$;

\smallskip
\noindent-\enskip\enskip if  $\vert j-j'\vert >1$, then $\beta_j$ and $\beta_{j'}$ are disjoint;

\smallskip
\noindent-\enskip\enskip if $\vert j-j'\vert =1$, then $\beta_j$ and $\beta_{j'}$ have a unique point of intersection and intersect transversely.

\medskip

We define a sequence $(\beta_j)_{-g\leq j\leq -1}$ in the other component of $M\setminus (\beta_{0-}\cup\beta_{0+})$ by setting 
$\beta_j=I^*(\beta_{-j})$.

\medskip
Here again we choose for every $j\in\{-g,\dots, -1,0-,0+,1\dots, g\}$, a closed annulus ${\cal A}_j^*$ containing $\beta_j$. We suppose that the annulus are disjoint except ${\cal A}^*_j$ and ${\cal A}^*_{j'}$ if $\vert j-j'\vert\leq 1$. For every $j\in\{-g,\dots, -1,0-,0+,1\dots, g\}$, we consider a Dehn twist $F_j^*$ supported on ${\cal A}_j^*$ and set $F_j={\cal E}([F_j^*])$. We define  $I={\cal E}(I^*)$ and $J={\cal E}(J^*)$, and  we know that  $I$ and $J$ are two commuting involutions. We know that $F_{0-}$ and $F_{0+}$ are conjugate by $J$ and  commute, that they commute with every $F_j$,  $j\in\{-g,\dots, -1,1,\dots,g\}$, and also with $I$. We know that $F_j$ and $F_{j'}$ commute if $\vert j-j'\vert>1$, that $F_j$ commutes with $J$ and is conjugate to $F_{-j}$ by $I$. The algebraic relation that should be true and that will give us a contradiction can be written
$$(F_1\circ \dots \circ F_g)^{g+1}=F_{0-}\circ F_{0+}=(F_{-1}\circ \dots \circ F_{-g})^{g+1}.$$

Like in the even case,  we know that for every $j\in\{-g,\dots, -1,1,\dots,g\}$, the domain of the canonical invariant essential cellular decomposition of  $F_j$ has a connected component $U_i$ which is a two punctured surface of genus $g-1$ and whose complement is an annular set ${\cal A}_j$,  intersection of a nested sequence of closed annuli homotopic to ${\cal A}_j^*$. We denote by ${\cal D}^j$ the decomposition restricted to $U_j$.  We have a similar situation for $F_{0-}$ and $F_{0+}$ and can define $U_{0-}$, $U_{0+}$, ${\cal A}_{0-}$,  ${\cal A}_{0+}$, ${\cal D}^{0-}$ and ${\cal D}^{0+}$. We can be more precise. The closure $Q$ of the complement of 
${\cal A}_{0-}^{*}\cup{\cal A}_{0+}^{*}$ has two connected components  that are compact two punctured surfaces of genus $(g-1)/2$. Let us consider the group generated by $F_{0-}$ and $F_{0+}$. By the remark stated at the end of Section 3, one can make prove, like it is done in Proposition 4.2, that the domain of the invariant canonical essential cellular decomposition of this group contains a surface with boundary isotopic to $Q$. This domain is included both in $U_{0-}$ and $U_{0+}$. One deduces that the sets ${\cal A}_{0-}$ and ${\cal A}_{0+}$ are disjoint and that $U_0=U_{0-}\cap U_{0+}$ has two connected components which are two punctured surfaces of genus $(g-1)/2$. Note that $I$ lets invariant  ${\cal A}_{0-}$ and ${\cal A}_{0+}$ and permutes the two connected components of $U_0$
and that $J$ permutes ${\cal A}_{0-}$ and  ${\cal A}_{0+}$ and lets invariant the two connected components of $U_0$. Let us recall what has been done in the even case.

\medskip
 \noindent{\bf i)}\enskip\enskip Lemma 4.3 is still valid. It can be applied to  ${\cal A}_{0-}$ and $F_{0-}$, and also to ${\cal A}_{0+}$ and $F_{0+}$. In particular, one can define  a set $ {\cal X}_{0,0-}\subset {\cal A}_{0-}$ and a set $ {\cal X}_{0,0+}\subset {\cal A}_{0+}$ as it was done in the even genus case.

\medskip
 \noindent{\bf ii)}\enskip\enskip  The sets $  {\cal X}_{0,0-}$ and $ {\cal X}_{0,0+}$ are non empty and permuted by $J$. The last point can be deduced from the fact that $J$ conjugates $F_{0-}$ to $F_{0+}$. To show the first point, observe that ${\rm Fix} (I)$, like ${\rm Fix} (I^*)$, contains four elements and is included in ${\cal A}_{0-}\cup {\cal A}_{0+}$. The set ${\rm Fix} (I)$ is invariant by $J$ because $I$ and $J$ commute. From the fact that  $J$ permutes ${\cal A}_{0-}$ and ${\cal A}_{0+}$, one deduces that there are two points of ${\rm Fix} (I)$ in ${\cal A}_{0-}$ and two points in ${\cal A}_{0+}$. One can conclude like in the even genus case.

\medskip
 \noindent{\bf iii)}\enskip\enskip  From the proof of Lemma 4.6, one deduces that for every $j\in \{-g,\dots,-1,1,\dots,g\}$, one has $ {\cal X}_{0,0-}\cup {\cal X}_{0,0+}\subset U_j$, and also that  $ {\cal X}_{0,0-}\subset U_{0+}$ and $ {\cal X}_{0,0+}\subset U_{0-}$.

\medskip
 \noindent{\bf iv)}\enskip\enskip  From the proof of Lemma 4.7,  one deduces that for every closed set $X\subset  {\cal X}_{0,0-}$ invariant by $F_{0-}$ and for every $j\in \{-g,\dots,-1,0+, 1,\dots g\}$, one has  ${\cal D}^j(X)\subset {\cal X}_{0,0-}$. One has a similar result changing $ {\cal X}_{0,0-}$ and $F_{0-}$ by  $ {\cal X}_{0,0+}$ and $F_{0+}$.

\medskip
 \noindent{\bf v)}\enskip\enskip   The arguments given in the proof of Lemma 4.8 and the fact that $J$ conjugates $F_{0-}$ to $F_{0+}$ tell us that if $X\subset  {\cal X}_{0,0-}$ is a closed connected set invariant by $F_{0-}$, then $X\cup J(X)$ separates the two connected components of $U_0$.

\medskip
 \noindent{\bf vi)}\enskip\enskip Let us consider the alphabet $\Xi^-=\{-g,\dots,-1,0+, 1,\dots g\}$ and a word $\overline j =(j_k)_{k\geq 0}\in (\Xi^{-})^{\N}$ in which every letter appears infinitely many often. One can define a non decreasing sequence $(X_k)_{k\geq 0}$ of closed sets invariant by $F_{0-}$ and included in $ {\cal X}_{0,0-}$,
by the induction formula
$$X_0={\rm Fix}(I)\cap {\cal A}_{0-}, \enskip X_{k+1}={\cal D}^{j_k}(X_k)\enskip \mathrm{if}\enskip k\geq 0.$$
The sequence $(X_k)_{k\geq 0}$ converges for the Hausdorff topology to $X_{\infty}=\overline{\bigcup_{k\geq 0} X_k}$. The set $X_{\infty}$ is invariant by $F_j$ for every $j\in\{-g,\dots,-1,0-,0+, 1,\dots g\}$. Moreover, there exists $k_*\geq 0$ such that $X_k$ has two connected components if $k<k_*$ and is connected if $k\geq k_*$.  We know that $X_k\cup J(X_k)$ separates the two connected components of $U_0$, as soon as $k\geq k_*$ and that it is the same for $X_{\infty}\cup J(X_{\infty})$.

\medskip
 \noindent{\bf vii)}\enskip\enskip We fix a connected component $V$ of $U_0$ and for every $k\geq K$, we write $V_k$ for the connected component of $M\setminus X_k\cup J(X_k)$ that contains $V$. It is a two punctured surface of genus $(g-1)/2$. The ``frontier'' of one of its ends is included in ${\cal A}_{0-}$. We consider its prime end compactification. We can define the corresponding rotation number $\rho^-_{\rm pe}(F_{0-}\vert_{V_k})$ and prove like in Lemma 4.9, that it is equal to $1/2+\Z$. Similarly we define the rotation number $\rho^-_{\rm pe}(F_{j_k}\vert_{V_k})$ and prove that it is equal to $0+\Z$.  Similarly, the connected component of  $M\setminus X_{\infty}\cup J(X_{\infty})$ that contains $V$ is a two punctured surface of genus $(g-1)/2$ and one of its ends has a frontier included in ${\cal A}_{0-}$. We define the rotation number  $\rho^-_{\rm pe}(F_{0-}\vert_{V_{\infty}})$ and prove like in Lemma 4.10, that it is equal to $1/2+\Z$. If $j\in \{-g,\dots,-1,0+, 1,\dots g\}$ we can define the rotation number $\rho^-_{\rm pe}(F_{j}\vert_{V_{\infty}})$ and show that it is equal to $0+\Z$.  We write $\Phi_j$, $j\in\{-g\dots,-1,0-,0+,1,\dots g\}$, for the homeomorphism defined on the circle of prime ends. 
 
\medskip
 \noindent{\bf viii)}\enskip\enskip From the fact that $F_{0-}$ and $F_{0+}$ commute, we deduce that $\Phi_{0-}$ and $\Phi_{0+}$ commute. One deduces that the rotation number of $\Phi_{0-}\circ \Phi_{0+}$ is the sum of the rotation number of $\Phi_{0-}$ and $\Phi_{0+}$, that means $1/2+\Z$. Therefore, $\Phi_{0-}\circ\Phi_{0+}$ has no fixed point. As we know that $\Phi_{0-}$ and $\Phi_{0+}$ commute with all the $\Phi_j$, $j\in\{-g\dots,-1,1,\dots g\}$, that every  $\Phi_j$, $j\in\{-g\dots,-1,1,\dots g\}$, has fixed points,  that $\Phi_j$ and $\Phi_{j'}$ commute if $\vert j-j'\vert >1$ and that
$$(\Phi_1\circ \dots \circ \Phi_g)^{g+1}=\Phi_{0-}\circ\Phi_{0+}=(\Phi_{-1}\circ \dots \circ \Phi_{-g})^{g+1}$$
here again, we can find a contradiction.

\bibliographystyle{alpha}
\renewcommand{\refname}{

  \end{document}